\documentclass{article}
\usepackage{amssymb}
\usepackage{amsmath,color}
\usepackage[hidelinks]{hyperref}

\usepackage[square,numbers]{natbib}         
\newtheorem{lem}{Lemma}[section]
\newtheorem{cor}[lem]{Corollary}
\newtheorem{teo}[lem]{Theorem}
\newtheorem{os}[lem]{Remark}

\newtheorem{prop}[lem]{Proposition}

\newcommand{\qed}{\thinspace\null\nobreak\hfill\hbox{\vbox{\kern-.2pt\hrule
 height.2pt depth.2pt\kern-.2pt\kern-.2pt \hbox to2.5mm{\kern-.2pt\vrule
 width.4pt \kern-.2pt\raise2.5mm\vbox to.2pt{}\lower0pt\vtop
 to.2pt{}\hfil\kern-.2pt \vrule
 width.4pt \kern-.2pt}\kern-.2pt\kern-.2pt\hrule height.2pt depth.2pt
 \kern-.2pt}}\par\medbreak}

\newcommand{\R}{\mathbb{R}}

\newcommand{\C}{\mathbb{C}}

\newcommand{\N}{\mathbb{N}}

\newcommand{\Rp}{\textrm{\emph{Re}\,}}

\newcommand{\eps}{\varepsilon}

\newcommand{\ds}{\displaystyle}

\date{}
\begin{document}

\title{Degenerate operators on  the half-line}
\author{G. Metafune \thanks{Dipartimento di Matematica e Fisica ``Ennio De Giorgi'', Universit\`a del Salento, C.P.193, 73100, Lecce, Italy.
-mail:  giorgio.metafune@unisalento.it}\qquad L. Negro \thanks{Dipartimento di Matematica e Fisica  ``Ennio De
Giorgi'', Universit\`a del Salento, C.P.193, 73100, Lecce, Italy. email: luigi.negro@unisalento.it} \qquad C. Spina \thanks{Dipartimento di Matematica e Fisica``Ennio De Giorgi'', Universit\`a del Salento, C.P.193, 73100, Lecce, Italy.
e-mail:  chiara.spina@unisalento.it}}

\maketitle
\begin{abstract}
\noindent 
We study elliptic and parabolic problems governed by the  singular elliptic   operators 
 \begin{equation*}
y^{\alpha}\left(D_{yy}+\frac{c}{y}D_y\right)-V(y),\qquad\alpha \in\R
\end{equation*}
in  $\R_+$, where  $V$ is  a potential having non-negative real part.
 
\par

\end{abstract}

\section{Introduction}
 In this paper we study solvability and regularity of elliptic and parabolic problems associated to the  degenerate   operators 
\begin{equation*} \label{defL}
 L =y^{\alpha}\left(D_{yy}+\frac{c}{y}D_y\right)-V \quad {\rm and}\quad D_t-  L
\end{equation*}
in the half-line $\R_+$.

Here $c,\alpha$ are real numbers and $V\in L^1_{loc}\left(\R^+,y^{c-\alpha}\right)$ is  a potential having non-negative real part.  The operator $B=D_{yy}+\frac{c}{y}D_y$ is a Bessel operator and satisfies the scaling property $$I_s^{-1} B I_s=s^{2-{\alpha}} B, \quad I_su(y)=u(sy).$$

\medskip
We study $L$  in  the weighted spaces $L^p_m:=L^p\left(\R^+,y^mdy\right)$, $m\in\R$,  and we characterize all $m$ such that  $L$ generates a $C_0$-semigroup. When  $V\geq 0$ we also prove that the generated semigroup is analytic   and we  show that it has maximal regularity, which means that both $D_t v$ and $L v$ have the same regularity as $(D_t - L) v$. In the case $V(y)=y^\alpha$ we finally characterize the domain of $L$.

We observe that the results already available for  $B$, see \cite[Section 3]{MNS-Caffarelli} and also  \cite{MNS-Max-Reg,MNS-Grad, MNS-Sharp,Negro-Spina-Asympt, met-negro-soba-spina} for  the $N$-d version of $B$, imply the corresponding ones for $y^\alpha B$ in $L^p_m$ by a change of variables, as described in Section \ref{Section Degenerate}.  The change of variables varies the underlying measure and explains why need the full scale of $L^p_m$ spaces.

More effort is needed to add the potential term.  We  consider first  $B-V$ in  $L^2(\R_+; y^c  dy)$. We use form methods to construct an analytic semigroup, we prove kernel bounds for complex times via Davies-Gaffney estimates  and provide a core. Then, with the methods of Section \ref{Section Degenerate}, we deduce similar results for $y^\alpha B-V$ in $L^2(\R_+; y^{c-\alpha}  dy)$.  Next   we prove that the semigroup can be extended to $L^p_m$ under sharp conditions on $p$ and $m$. 
Finally  we prove that for every $\epsilon >0$  the family of operators  
\begin{align*}
	\left\{e^{z(y^\alpha B-V)}:\; z\in \Sigma_{\frac \pi 2-\epsilon},\; 0\leq V\in L^1_{loc}\left(\R^+,y^{c-\alpha}\right)\right\}
\end{align*}
is $\mathcal R$-bounded in $L^p_m$, which implies the maximal regularity  of the semigroup when $V\geq 0$.

As a motivation for our investigation, we point out that, in the special case $V(y)=y^\alpha$,  all the results above   play   a crucial role in \cite{MNS-CompleteDegenerate} in  the investigation  of the degenerate operators
	\begin{align*}
		\mathcal L =y^{\alpha_1}\Delta_{x} +y^{\alpha_2}\left(D_{yy}+\frac{c}{y}D_y  -\frac{b}{y^2}\right).
	\end{align*}
Let us suppose, for simplicity,  $b=0$, $\alpha_1=\alpha_2:=\alpha$.	 Assuming that $y^\alpha(\Delta_x u+B_yu)=f$ and taking the Fourier transform ${\cal F} u$ or $\hat u$  with respect to $x$ we obtain  $y^\alpha|\xi|^2 \hat u(\xi,y)=-y^\alpha |\xi|^2 (y^\alpha |\xi|^2-y^\alpha B_y)^{-1}\hat f (\xi,y)$. 
	Therefore
\begin{align*}
		y^\alpha\Delta_x  \mathcal  L^{-1}&={\cal F}^{-1}\left (y^\alpha|\xi|^2 (y^\alpha|\xi|^2-y^\alpha B_y)^{-1}) \right) {\cal F}
\end{align*}
	and the boundedness of  $y^\alpha\Delta_x  \mathcal L^{-1}$ is equivalent to that of the multiplier
	\begin{align*}
		\xi \in \R^N \to y^\alpha|\xi|^2(y^\alpha|\xi|^2-y^\alpha By)^{-1}.
	\end{align*}
For this reason we prove in Section \ref{case V yalpha} that certain multiplier associated to $y^\alpha B- V$  satisfy a vector valued Mikhlin  theorem. These results rely on  square function estimates which we deduce from kernel bounds
and the following equality, which allows to treat $\lambda$  or $|\xi|^2$ as spectral parameters simultaneously
\begin{align*}
	\left(\lambda-y^\alpha B+|\xi|^2 y^\alpha\right)^{-1}=\left(|\xi|^2- B+\frac{\lambda}{y^\alpha}\right)^{-1}\frac{1}{y^{\alpha}}.
\end{align*}


We restrict ourselves to $\alpha<2$ and  consider $y^\alpha B$ with Neumann boundary condition at $0$, namely  $\lim_{y\to 0} y^cD_yu(y)=0$. This is equivalent to require   $y^{\alpha-1}D_{y}u\in L^p_m$, see \cite[Proposition 5.11]{MNS-Sobolev}. The restriction $\alpha<2$ is not really essential
since one can deduce from it the case $\alpha>2$, which requires a  boundary condition at $\infty$, using the change of variables described in Section \ref{Section Degenerate}.

Besides this, our strategy  can be  easily adapted to different boundary conditions and to more general operators $y^\alpha\left(D_{yy}+\frac{c}{y}D_y-\frac{b}{y^2}\right)-V$. 
We do this  (in much more generality) in \cite[Sections 7,8]{MNS-CompleteDegenerate}. 

\medskip

The paper is organized as  follows. In Section 2 we briefly recall the harmonic analysis background needed in the paper,  as  square function estimates, $\mathcal R$-boundedness  and a vector valued multiplier theorem. 

In Section 3, we exploit an elementary change of variables, in a functional analytic setting, to reduce our operators to the simpler case where $\alpha=0$.

Section 4 is devoted to the study of the Bessel operator $y^\alpha B$. In Sections 5, 6 and 7 we perturb the Bessel operator  by adding the potential $V$ and we prove real and complex  kernel estimates, generation results and maximal regularity for $y^\alpha B-V$. Finally in Section 8 we treat the case $V(y)=y^\alpha$ and characterize the domain of  $y^\alpha B-y^{\alpha}$.

\bigskip
\noindent\textbf{Notation.} For $m \in \R$ we consider the measure $y^m dy $ in $\R_+$ and  we write $L^p_m$ for  $L^p(\R_+, y^m  dy)$. Similarly $W^{k,p}_m=\{u \in L^p_m: \partial^\alpha u \in  L^p_m \quad |\alpha| \le k\}$. When we write $V \in L^q_{loc}(\R^+, y^m\, dy)$, we mean that $V\in  L^q([0,b], y^m\, dy)$ for every $b <\infty$. 

We use $\C^+=\{ \lambda \in \C: \Rp \lambda >0 \}$ and, for $|\theta| \leq \pi$, we denote by  $\Sigma_{\theta}$  the open sector $\{\lambda \in \C: \lambda \neq 0, \ |Arg (\lambda)| <\theta\}$.
\bigskip

\section{Harmonic analysis and maximal regularity}
The study of maximal regularity of parabolic problems of the form $u_t=Au+f, u(0)=0$, where $A$ is the generator of an analytic semigroup on a Banach space $X$, consists in proving estimates like
$$
\|u_t\|_p+\|Au\|_p \le \|f\|_p
$$
where the $L^p$ norm is that of $L^p([0, T[;X)$. This  can be interpreted as closedness of $D_t-A$ on the intersection of the respective domains or, equivalently, boundedness of the operator $A(D_t-A)^{-1}$ in $L^p([0, T[;X)$.

Nowadays this strategy is well established and relies on Mikhlin vector-valued multiplier theorems.
Let us state the relevant definitions and main results we need, referring the reader to \cite{DenkHieberPruss}, \cite{Pruss-Simonett} or \cite{KunstWeis}.

Let ${\cal S}$ be a subset of $B(X)$, the space of all bounded linear operators on a Banach space $X$. ${\cal S}$ is $\mathcal R$-bounded if there is a constant $C$ such that
$$
\|\sum_i \eps_i S_i x_i\|_{L^p(\Omega; X)} \le C\|\sum_i \eps_i  x_i\|_{L^p(\Omega; X)} 
$$
for every finite sum as above, where $(x_i ) \subset X, (S_i) \subset {\cal S}$ and $\eps_i:\Omega \to \{-1,1\}$ are independent and symmetric random variables on a probability space $\Omega$. The smallest constant $C$ for which the above definition holds is the $\mathcal R$-bound of $\mathcal S$,  denoted by $\mathcal R(\mathcal S)$.
It is well-known that this definition does not  depend on $1 \le p<\infty$ (however, the constant $\mathcal R(\mathcal S)$ does) and that $\mathcal R$-boundedness is equivalent to boundedness when $X$ is an Hilbert space.
When $X$ is an $L^p$ space (with respect to any $\sigma$-finite measure), testing  $\mathcal R$-boundedness is equivalent to proving square functions estimates, see \cite[Remark 2.9 ]{KunstWeis}.

\begin{prop}\label{Square funct R-bound} Let ${\cal S} \subset B(L^p(\Sigma))$, $1<p<\infty$. Then ${\cal S}$ is $\mathcal R$-bounded if and only if there is a constant $C>0$ such that for every finite family $(f_i)\in L^p(\Sigma), (S_i) \in {\cal S}$
$$
\left\|\left (\sum_i |S_if_i|^2\right )^{\frac{1}{2}}\right\|_{L^p(\Sigma)} \le C\left\|\left (\sum_i |f_i|^2\right)^{\frac{1}{2}}\right\|_{L^p(\Sigma)}.
$$
\end{prop}
The best constant $C$ for which the above square functions estimates hold satisfies $\kappa^{-1} C \leq \mathcal R(\mathcal S) \leq \kappa C$ for a suitable $\kappa>0$ (depending only on $p$). The proposition above $\mathcal R$-boundedness follows from domination.
\begin{cor} \label{domination}
Let  ${\cal S}, {\cal T} \subset B(L^p(\Sigma))$, $1<p<\infty$ and assume that $\cal T$ is $\mathcal R$ bounded and that for every $S \in \cal S$ there exists $T \in \cal T$ such that $|Sf| \leq |Tf|$ pointwise, for every $f \in L^p(\Sigma)$. Then ${\cal S}$ is $\mathcal R$-bounded.
\end{cor}

Let $(A, D(A))$ be a sectorial operator in a Banach space $X$; this means that $\rho (-A) \supset \Sigma_{\pi-\phi}$ for some $\phi <\pi$ and that $\lambda (\lambda+A)^{-1}$ is bounded in $\Sigma_{\pi-\phi}$. The infimum of all such $\phi$ is called the spectral angle of $A$ and denoted by $\phi_A$. Note that $-A$ generates an analytic semigroup if and only if $\phi_A<\pi/2$. The definition of $\mathcal  R$-sectorial operator is similar, substituting boundedness of $\lambda(\lambda+A)^{-1}$ with $\mathcal R$-boundedness in $\Sigma_{\pi-\phi}$. As above one denotes by $\phi^R_A$ the infimum of all $\phi$ for which this happens; since $\mathcal R$-boundedness implies boundedness, we have $\phi_A \le \phi^R_A$.

\medskip

The $\mathcal R$-boundedness of the resolvent characterizes the regularity of the associated inhomogeneous parabolic problem, as we explain now.

An analytic semigroup $(e^{-tA})_{t \ge0}$ on a Banach space $X$ with generator $-A$ has
{\it maximal regularity of type $L^q$} ($1<q<\infty$)
if for each $f\in L^q([0,T];X)$ the function
$t\mapsto u(t)=\int_0^te^{-(t-s)A})f(s)\,ds$ belongs to
$W^{1,q}([0,T];X)\cap L^q([0,T];D(B))$.
This means that the mild solution of the evolution equation
$$u'(t)+Au(t)=f(t), \quad t>0, \qquad u(0)=0,$$
is in fact a strong solution and has the best regularity one can expect.
It is known that this property does not depend on $1<q<\infty$ and $T>0$.
A characterization of maximal regularity is available in UMD Banach spaces, through the $\mathcal  R$-boundedness of the resolvent in a suitable sector $\omega+\Sigma_{\phi}$, with $\omega \in \R$ and $\phi>\pi/2$ or, equivalently, of the scaled semigroup $e^{-(A+\omega')t}$ in a sector around the positive axis. In the case of $L^p$ spaces it can be restated in the following form,  see \cite[Theorem 1.11]{KunstWeis}

\begin{teo}\label{MR} Let $(e^{-tA})_{t \ge0}$ be a bounded analytic semigroup in $L^p(\Sigma)$, $1<p<\infty$,  with generator $-A$. Then $T(\cdot)$ has maximal regularity of type $L^q$  if and only if the set $\{\lambda(\lambda+A)^{-1}, \lambda \in  \Sigma_{\pi/2+\phi} \}$ is $\mathcal R$- bounded for some $\phi>0$. In an equivalent way, if and only if 
there are constants $0<\phi<\pi/2 $, $C>0$ such that for every finite sequence $(\lambda_i) \subset \Sigma_{\pi/2+\phi}$, $(f_i) \subset  L^p$
$$
\left\|\left (\sum_i |\lambda_i (\lambda_i+A)^{-1}f_i|^2\right )^{\frac{1}{2}}\right\|_{L^p(\Sigma)} \le C\left\|\left (\sum_i |f_i|^2\right)^{\frac{1}{2}}\right\|_{L^p(\Sigma)}
$$
or, equivalently, there are constants $0<\phi'<\pi/2 $, $C'>0$ such that  for every finite sequence
$(z_i) \subset \Sigma_{\phi'}$, $(f_i) \subset  L^p$
$$
\left\|\left (\sum_i |e^{-z_i A}f_i|^2\right )^{\frac{1}{2}}\right\|_{L^p(\Sigma)} \le C'\left\|\left (\sum_i |f_i|^2\right)^{\frac{1}{2}}\right\|_{L^p(\Sigma)}.
$$
\end{teo}
\medskip
Finally we state  a version of the operator-valued Mikhlin multiplier theorem in the N-dimensional case, see \cite[Theorem 3.25]{DenkHieberPruss} or \cite[Theorem 4.6]{KunstWeis}.
\begin{teo}   \label{mikhlin}
Let $1<p<\infty$, $M\in C^N(\R^N\setminus \{0\}; B(L^p(\Sigma))$ be such that  the set
$$\left \{|\xi|^{|\alpha|}D^\alpha_\xi M(\xi): \xi\in \R^{N}\setminus\{0\}, \ |\alpha | \leq N \right \}$$
is $\mathcal{R}$-bounded.
Then the operator $T_M={\cal F}^{-1}M {\cal F}$ is bounded in $L^p(\R^N, L^p(\Sigma))$, where $\cal{F}$ denotes the Fourier transform.
\end{teo}

We end this section with the following lemma on radially symmetric multipliers.
\begin{lem}   \label{mikhlin Comp}
	Let $1<p<\infty$, $m\in C^N(\R_+; B(L^p(\Sigma))$ be such that  the set
	$$\left \{s^{k} m^{(k)}(s): \ s\in \R_+, \ k\leq N \right \}$$
	is $\mathcal{R}$-bounded. For $a\in\R$ let 
	$M(\xi)=m\left(|\xi|^a\right)$.
Then  $M\in C^N(\R^N\setminus \{0\}; B(L^p(\Sigma))$ and  
	$$\left \{|\xi|^{|\alpha|}D^\alpha_\xi M(\xi): \xi\in \R^{N}\setminus\{0\}, \ |\alpha | \leq N \right \}$$
	is $\mathcal{R}$-bounded and  
$$\mathcal R\left \{|\xi|^{|\alpha|}D^\alpha_\xi M(\xi): \xi\in \R^{N}\setminus\{0\}, \ |\alpha | \leq N \right \}\leq C(N)\mathcal R\left \{s^{k} m^{(k)}(s): \ s\in \R_+, \ k\leq N \right \}.$$
\end{lem}
{\sc{Proof.}} Let us observe preliminarily that for any multi-index $\alpha$ with $0<|\alpha|\leq N$ one has
\begin{align}\label{eq mik comp}
	D^{\alpha}_\xi M(\xi)
		&=
		\ds \sum_{i=1}^{|\alpha|} h_{i,\alpha}(\xi)m^{(i)}\left(|\xi|^a\right)
\end{align}
where $h_{i,\alpha} \in C^\infty (\R^N\setminus \{0\})$ are homogeneous functions of degree $ia-|\alpha|$.  Obviously $\eqref{eq mik comp}$ is valid for $|\alpha|=1$ since $\nabla M(\xi)=a \ m'(|\xi|^a)|\xi|^{a-2}\xi$ and follows by induction, since the derivatives of  $h_{i,\alpha}$ are homogeneous of degree $ia-|\alpha|-1$. 

The proof of the lemma now follows by Corollary  \eqref{domination} since from \eqref{eq mik comp} one has  for $f \in L^p(\Sigma)$
\begin{align*}
|\xi|^{|\alpha|}|D^\alpha_\xi M(\xi)f|&\leq \ds |\xi|^{|\alpha|}\sum_{i=1}^{|\alpha|} |h_{i,\alpha}(\xi)| |m^{(i)}\left(|\xi|^a\right)f|\leq C\ds\sum_{i=1}^{|\alpha|} |\xi|^{ia} |m^{(i)}\left(|\xi|^a\right)f|.
\end{align*}
\qed

\section{Degenerate operators and similarity transformations }\label{Section Degenerate}
We investigate when the  operators 
$$
B=D_{yy}+\frac{c}{y}D_y, \qquad y^\alpha B= y^\alpha\left(D_{yy}+\frac{c}{y}D_y\right)
$$
can be transformed one into the other by means of  change of variables. Here $\alpha,c$ are unrestricted real coefficients. 

For  $\beta \in\R$, $\beta\neq -1$ let 
\begin{align}\label{Gen Kelvin def}
T_{\beta\,}u(y)&:=|\beta+1|^{\frac 1 p}u(y^{\beta+1}),\quad y\in\R_+.
\end{align}
Observe that
$$ T_{\beta\,}^{-1}=T_{-\frac{\beta}{\beta+1}\,}.$$

\begin{prop}\label{Isometry action der} Let $1\leq p\leq \infty$, $k,\beta \in\R$, $\beta\neq -1$. The following properties hold.
	\begin{itemize}
		\item[(i)] For every $m\in\R$,  $T_{\beta\,}$ maps isometrically  $L^p_{\tilde m}$ onto $L^p_m$  where 
		$$ \tilde m=\frac{m-\beta}{\beta+1}.$$
		\item[(ii)] For every  $u\in W^{2,1}_{loc}\left(\R_+\right)$ one has
\begin{itemize}
	\item[1.] $y^\alpha T_{\beta\,}u=T_{\beta\,}(y^{\frac{\alpha}{\beta+1}}u)$, for any $\alpha\in\R$;\medskip
	\item[2.]  $D_y T_{\beta\,}u=T_{\beta\,}\left((\beta+1)y^{\frac{\beta}{\beta+1}}D_yu\right)$,
	\\[1ex] $D_{yy} (T_{\beta\,} u)=T_{\beta\,}\Big((\beta+1)^2y^{\frac{2\beta}{\beta+1}}D_{yy}u+(\beta+1)\beta y^{\frac{\beta-1}{\beta+1}}D_y u\Big)$.\medskip
\end{itemize}
	\end{itemize}
\end{prop}{\sc{Proof.}} The proof of (i) follows after observing the Jacobian of $y\mapsto y^{\beta+1}$ is $|1+\beta|y^{\beta}$. 
Then we compute  
\begin{align*}
D_y T_{\beta\,}u(y)=&|\beta+1|^{\frac 1 p}\left((\beta+1)y^\beta D_y u(y^{\beta+1})\right)=T_{\beta\,}\left((\beta+1)y^{\frac{\beta}{\beta+1}}D_yu\right)
\end{align*}
and similarly
\begin{align*}
D_{yy} T_{\beta\,} u(y)=&T_{\beta\,}\Big((\beta+1)^2y^{\frac{2\beta}{\beta+1}}D_{yy}u+(\beta+1)\beta y^{\frac{\beta-1}{\beta+1}}D_y u\Big).
\end{align*}
\qed
\begin{prop}\label{Isometry action}  Let 
	$T_{\beta\,}$ be the isometry above defined.  The following properties holds.
	
		For every  $u\in W^{2,1}_{loc}\left(\R_+\right)$ one has
		$$
		T_{\beta\,}^{-1} \Big(y^{\alpha} B\Big)T_{\beta\,}u=\Big((\beta+1)^2y^{\frac{\alpha+2\beta}{\beta+1}}\tilde{ B}\Big) u
		$$
		where $\tilde { B}$ is the operator defined as in \eqref{defL} with parameter $c$ replaced, respectively, by
		$$
			\tilde c=\frac{c+\beta\left(c+1+\beta\right)}{(\beta+1)^2}.
		$$
\end{prop}
{\sc{Proof.}} Using Proposition \ref{Isometry action der} we can compute
\begin{align*}
	{B}T_{\beta\,}u(y)&=T_{\beta\,}\Big[(\beta+1)^2y^{\frac{2\beta}{\beta+1}}D_{yy}u+(\beta+1)\beta y^{\frac{\beta-1}{\beta+1}}D_y u+c(\beta+1)y^{\frac{\beta-1}{\beta+1}}D_yu-by^{-\frac{2}{\beta+1}}u \Big]\\[1ex]
	&=T_{\beta\,}\Bigg[y^{\frac{2\beta}{\beta+1}}
	\Bigg(
	(\beta+1)^2D_{yy} u+\frac{(\beta+1)\left(\beta+c\right)}{y}D_y u-b\frac u{y^2}
	\Bigg)
	\Bigg]=T_{\beta\,}\left(y^{\frac{2\beta}{\beta+1}}\tilde { B} u\right)
\end{align*}
which implies
\begin{align*}
	T_{\beta\,}^{-1}\left(y^{\alpha}{B}\right)T_{\beta\,}u&=y^{\frac{\alpha+2\beta}{\beta+1}}\tilde { B} u.
\end{align*}
\qed

\section{The Bessel operator $y^\alpha B^n$}\label{sec y^aB}
In this section we consider for $\alpha < 2$,  $c \in\R$  the  operator 
\begin{equation*} 
y^\alpha B =y^\alpha\left(D_{yy}+\frac{c}{y}D_y  \right)
\end{equation*}
in the space $L^p_m$ under Neumann boundary conditions. 

According to Proposition \ref{Isometry action},  we use the isometry 
$$T_{-\frac \alpha 2}:L^p_{\tilde m}\to L^p_m\,\quad T_{-\frac \alpha 2}u(y)=\left|1-\frac{\alpha}{2}\right |^{\frac 1p}u(y^{1-\frac{\alpha}{2}}),$$
$\tilde m=\frac{m+\frac\alpha 2}{1-\frac\alpha 2}$, 
under which $y^\alpha B$ becomes  isometrically  equivalent to  
$
T_{-\frac \alpha 2}^{-1} \Big(y^{\alpha}B\Big)T_{-\frac \alpha 2}=\left(1-\frac\alpha 2\right)^2\tilde B
$ where $\tilde B=D_{yy}+\frac{\tilde c}{y}D_y $ and  
$
 \tilde c=\frac{c-\frac\alpha 2}{1-\frac\alpha 2}$.

 All the results for $y^\alpha B$ in $L^p_m$ are then  immediate consequence of those of $\tilde{B}$ in $L^p_{\tilde m}$ already proved in \cite[Section 3]{MNS-Caffarelli} (see also \cite{MNS-Max-Reg,MNS-Grad, MNS-Sharp,Negro-Spina-Asympt} for some analogous results for the $N$-d version of $\tilde B$).

  If $1<p<\infty$, we define 
\begin{equation*} \label{w2p}
W^{2,p}_{\mathcal N}(\alpha,m)=\left\{u\in W^{2,p}_{loc}(\R_+):\ u,\    y^{\alpha}D_{yy}u,\ y^{\frac{\alpha}{2}}D_{y}u,\ y^{\alpha-1}D_{y}u\in L^p_m\right\}
\end{equation*}
and refer to \cite{MNS-Sobolev} where these spaces are studied in detail in $\R^{N+1}_+$. 
The Neumann boundary condition is enclosed in the requirement $y^{\alpha-1}D_{y}u\in L^p_m$. This last is redundant when $(m+1)/p>1-\alpha$ and equivalent to $D_yu(y) \to 0$ as $y \to 0$, when $(m+1)/p<1-\alpha$, see \cite[Proposition 4.3]{MNS-Sobolev}.
Consequently we write $y^\alpha B^n$ or, more pedantically $y^\alpha B^n_{m,p}$ if necessary, for its realization in $L^p_m$.

 \begin{os}
 The restriction $\alpha<2$ is not really essential
 since one can deduce from it the case $\alpha>2$, which requires boundary condition at $\infty$, using the change of variables described in Section \ref{Section Degenerate} or directly from the equality $
 T_{-\frac \alpha 2}^{-1} \Big(y^{\alpha}B\Big)T_{-\frac \alpha 2}=\left(1-\frac\alpha 2\right)^2\tilde B
 $ which is valid for any $\alpha\neq 2$. However, here and in what follows, we keep to it in order  to simplify the exposition. 
\end{os}

\begin{teo}\label{Neumnan alpha m 2}
If  $0<\frac{m+1}p<c+1-\alpha$, then 
$y^\alpha B^n$ endowed with domain $W^{2,p}_{\mathcal N}(\alpha,m)$

generates a bounded positive analytic semigroup of angle $\pi/2$ on $L^p\left(\R_+,y^mdy\right)$. 

\end{teo}
{\sc Proof.} We use the identity  $T_{-\frac \alpha 2}^{-1} \Big(y^{\alpha}B^n\Big)T_{-\frac \alpha 2}=\left(1-\frac\alpha 2\right)^2\tilde{B^n}$ and apply \cite[Proposition 3.3]{MNS-Caffarelli} in $L^p_{\tilde m}$. Note that   $D(y^\alpha B_{m,p}^n)=T_{-\frac \alpha 2}D( \tilde {B}_{\tilde m,p}^n)$ which means
\begin{align*}
u\in D(y^\alpha B_{m,p}^n)\qquad\Longleftrightarrow\qquad v(y):=u(y^{\frac 2{2-\alpha}})\in D( \tilde{B}_{\tilde m,p}^n).
\end{align*}
\qed


Under the hypothesis of Theorem 	\ref{Neumnan alpha m 2}, the domain of $y^{\alpha}B^n$ consists of all functions in the maximal domain  satisfying a  Neumann  condition at $0$, see \cite[Proposition 4.6, 4.7]{MNS-Sobolev}, that is  
$$ D(y^\alpha B_{m,p}^n)=\left\{u \in  W^{2,p}_{loc}(\R_+): u,\  y^{\alpha}Bu \in L^p_m\text{\;\;and\;\;}\lim_{y\to 0}y^c D_yu=0\right\}$$
(the condition $\lim_{y\to 0}y^c D_yu=0$ can be deleted in the range $0<\frac{m+1}p\leq c-1$). When $c\geq 1$ the domain can also be described  involving  a Dirichlet, rather than Neumann, boundary condition
\begin{align*}
	 D(y^\alpha B_{m,p}^n)=&\left\{u \in  W^{2,p}_{loc}(\R_+): u,\ y^{\alpha}Bu \in L^p_m\text{\;and\;}\lim_{y\to 0}y^{c-1} u=0\right\},\qquad \text{\hspace{1ex}if $c>1$};\\[1.5ex]
		 D(y^\alpha B_{m,p}^n)=&\left\{u \in  W^{2,p}_{loc}(\R_+): u,\ y^{\alpha}Bu \in L^p_m\text{\;and\;}\lim_{y\to 0} u\in \C\right\},\qquad \text{\hspace{5.5ex}if $c=1$}.
\end{align*}



We close this section by describing a core  which does not depend on $\alpha,m,p$ and on the coefficients of the operator.
\begin{prop} \label{core}
If $0<\frac{m+1}{p}<c+1-\alpha$, then  a core for $y^{\alpha}B^n$ is 
$$\mathcal {D}= \left\{u \in C_c^\infty ([0,\infty)): u \ {\rm constant \ in \ a \ neighborhood\  of }\  0  \right \}.$$
\end{prop}
{\sc Proof.} The proof immediately follows by observing that, by \cite[Proposition 5.4]{MNS-Caffarelli}, $\mathcal {D}$ is a core when $\alpha=0$, that is for $\tilde{B}^n_{\tilde {m},p}$,  and the isometry $T_{-\frac{\alpha}2}$ leaves invariant $\mathcal{D}$ since  $\alpha<2$.
\qed

	\begin{os}\label{Core compact}
		We point out that, by  
		the proof of \cite[Proposition 5.4]{MNS-Caffarelli} or by \cite[Remark 4.14]{MNS-Sobolev},  it follows that  if $u\in D(y^{\alpha}B^n_{m,p})$ has  support in $[0,b]$, then there exists a sequence $\left(u_n\right)_{n\in\N}\in\mathcal D$  such that $ \mbox{supp }u_n\subseteq [0,b]$ and $u_n\to u$ in $D(y^{\alpha}B^n_{m,p})$.
	\end{os}

\section{The operator $B^n-V$}
We start our investigation by adding a potential $0 \leq V\in L^1_{loc}\left(\R^+,y^c\, dy\right)$ to $B^n$. Here we  prove kernel bounds and construct a core.
\subsection{Kernel bounds}
For $c+1>0$ and  $0 \leq V\in L^1_{loc}\left(\R^+,y^c\, dy\right)$ we  prove  upper bounds for the heat kernel of $ B^n-V$, following the method used in \cite[Sections 3,4]{met-calv-negro-spina}.

Setting $H^1_c=\{u \in L^2_c, u' \in L^2_c\}$,  we recall that from \cite[Section 2]{MNS-Caffarelli}  the operator $B^n_{c,2}$ is associated to the non-negative, symmetric and closed form  in $L^2_c$
\begin{align*}
\mathfrak{a}(u,v)
&:=
\int_0^\infty D_y u  D_y \overline{v}\,y^c dy,\qquad  D(\mathfrak{a})=H^1_c.
\end{align*}
We  consider the perturbed  form $\mathfrak{ a}_V$ in $L^2_{c}$ defined by
\begin{align}\label{form v}
\nonumber&\mathfrak{a}_V(u,v)=\mathfrak{a}(u,v)+ \left\langle Vu,v\right\rangle_{L^2_c}=
\int_{\R_+}
\left(D_y u D_y \overline{v}+Vu\overline{v}\right)y^c\,dy\\[1ex]
&D(\mathfrak{a}_V)=D(\mathfrak{a})\cap L^{2}\left(\R^+,Vy^c\, dy\right)
\end{align}
and define $B^n-V$ in $L^2_c$ as the operator  associated to the form $\mathfrak{ a}_V$
\begin{align*}
D(B^n-V)&=\{u \in D(\mathfrak{a}_V): \exists  f \in L^2_c \ {\rm such\ that}\  \mathfrak{a}_V(u,v)=\int_0^\infty f \overline{v}y^c\, dy\ {\rm for\ every}\ v\in D(\mathfrak{a}_V)\},\\
  B^nu-Vu&=-f.
\end{align*}
The positivity of $V$ implies that  the norm induced by the form $\mathfrak{a}_V$ is stronger then the one induced by $\mathfrak{a}$: as an  immediate consequence one deduces that $\mathfrak{a}_V$ is closed. 
By standard theory on sesquilinear forms we have  the following result.

\begin{prop}\label{kernel Vreal}
If $c+1>0$, $0 \leq V\in L^1_{loc}\left(\R^+,y^c\, dy\right)$,  then $\mathfrak{a}_V$ 
is a non-negative, symmetric and closed form in $L^2_{c}$. Its associated operator $-B^n+V$ is  non-negative and self-adjoint and  $B^n-V$ generates a contractive analytic semigroup $\left\{e^{z(B^n-V)}:\ z\in\C_+\right\}$  in $L^2_{c}$. Moreover:
\begin{itemize}
  \item[(i)] The semigroup $\left(e^{t(B^n-V)}\right)_{t\geq 0}$ is sub-markovian (i.e. it is positive and $L^\infty$-contractive)  and it is dominated by $e^{tB^n}$, that is 
\begin{align*}
 |e^{t(B^n-V)}f|\leq e^{tB^n}|f|,\qquad t>0,\quad  f\in L^2_c.
 \end{align*}
   \item[(ii)] $\left(e^{t(B^n-V)}\right)_{t\geq 0}$ is a semigroup of integral operators and its heat kernel $p_V$, taken with respect the measure $\rho^cd\rho$, satisfies
   \begin{align*}
0\leq p_{V}(t,y,\rho)\leq C t^{-\frac{1}{2}}  \rho^{-c}\left (\frac{\rho}{t^{\frac{1}{2}}}\wedge 1 \right)^{c}
\exp\left(-\frac{|y-\rho|^2}{\kappa t}\right).
    \end{align*}     
  \end{itemize}  
\end{prop}
{\sc Proof.} The first claim follows from the property of $\mathfrak{a}_V$.  $e^{t(B^n-V)}$ is sub-markovian   from  \cite[Corollary 2.17]{Ouhabaz}.
The domination property follows  from \cite[Corollary 2.21]{Ouhabaz}.
(ii) is a consequence of  \cite[Proposition 1.9]{ArendtBukhalov} since 
$e^{t(B^n-V)}$ is dominated by the positive integral operator $e^{tB^n}$ whose kernel satisfies the stated estimate, see  \cite[Proposition 2.8]{MNS-Caffarelli} where, however, the  kernel is written with respect to the Lebesgue measure.\qed

To  extend the above heat kernel estimates to the half plane $\C_+$ we need the following lemma.

\begin{lem}\label{Misura palle}
Let $c+1>0$ and for  $y_0,r>0$
\begin{align*}
Q_c (y_0,r):=\int_{[y_0,y_0+r]} y^{c} dy.
\end{align*}
 Then one has 
\begin{align*}
Q_c(y_0,r)\simeq r^{c+1}\left(\frac{y_0}{r}\right)^{c}\left(\frac{y_0}{r}\wedge 1\right)^{-c},\qquad r,y_0>0 .
\end{align*}
In particular the function $Q_c$ satisfies, for some  constants $C\geq 1$, the doubling condition
\begin{align*}
\frac{Q_c(y_0,s)}{Q_c(y_0,r)}\leq C \left(\frac{s}{r}\right)^{1\vee (c+1)},\qquad \forall y_0>0,\quad 0<r< s.
\end{align*}
\end{lem}
{\sc{Proof.}}  A scaling argument immediately yields $Q_c(y_0,r)=r^{c+1}Q_c\left(\frac{y_0}r,1\right)$ and we may therefore assume  $r=1$. The local integrability of $y^c$ implies that  $Q_c (y_0,1)$ is continuous as a function of $y_0$ and moreover $Q_c (y_0,1)\to \int_{(0,1)} y^{c} dy>0$ as $y_0\to 0$. Therefore   if  $y_0 \leq 1$ then 
\begin{align*}
Q_c (y_0,1)\simeq 1.
 \end{align*}
On the other hand if $y_0>1$ then $y\simeq y_0$ for any $y\in (y_0,y_0+1)$ which implies
\begin{align*}
Q_c (y_0,1)=\int_{(y_0,y_0+1)} y^{c} dy\simeq y_0^{c}.
 \end{align*}
The last two inequalities yields $Q_c (y_0,1)\simeq \left(y_0\right)^{c}\left(y_0\wedge 1\right)^{-c}$. The doubling condition follows from the previous estimates and the fact that for $0<r<s$ one has 
\begin{align*}
\frac{Q_c(y_0,s)}{Q_c(y_0,r)}\leq C \begin{cases}
\left(\frac{s}{r}\right)^{c+1},\qquad &\text{if}\quad \frac{y_0}s\leq \frac{y_0}r\leq 1;\\[1ex]
\frac{s}{r}(\frac{s}{y_0})^c,\qquad &\text{if}\quad \frac{y_0}s\leq 1< \frac{y_0}r;\\[1ex]
\frac{s}{r},\qquad &\text{if}\quad 1\leq \frac{y_0}s\leq \frac{y_0}r
\end{cases}
\end{align*}
(note that in the range $\frac{y_0}s\leq 1< \frac{y_0}r$ one has $(\frac{s}{y_0})^c\leq 1 $ if $c<0$ and $(\frac{s}{y_0})^c\leq \left(\frac{s}{r}\right)^{c}$ if $c\geq 0$).\qed

\begin{prop}\label{kernel V}
Let $c+1>0$, $0 \leq V\in L^1_{loc}\left(\R^+,y^c\, dy\right)$. The semigroup $\left\{e^{z(B^n-V)}:\ z\in\C_+\right\}$ consists of integral operators
\[
e^{z(B^n-V)}f(y)=
\int_{0}^\infty p_{V}(z,y,\rho)f(\rho)\,\rho^c d\rho,\quad f\in L^2_c,\quad y>0.
\]
Furthermore for  every $\epsilon>0$ there exist $k_\epsilon,C_\epsilon>0$ such that, for every  $z\in\Sigma_{\frac\pi 2-\epsilon}$ and $y,\rho>0$,
\begin{align*}
|p_V(z,x,y)|\leq
C_\epsilon |z|^{-\frac{1}{2}}  \rho^{-c}\left (\frac{\rho}{|z|^{\frac{1}{2}}}\wedge 1 \right)^{c}
\exp\left(-\frac{|y-\rho|^2}{\kappa_\epsilon |z|}\right).
\end{align*}
\end{prop}
{\sc Proof.}  Using the previous lemma  we rewrite Proposition \ref{kernel Vreal} (ii) as  
\begin{align*}
0\leq p_V(t,y,\rho)\leq C \frac{1}{Q_c(\rho,\sqrt t)}\exp\left(-\frac{|y-\rho|^2}{\kappa t}\right).
\end{align*}
Furthermore by \cite[Theorem 3.3]{Coulhon-Sikora},  $e^{t(B^n-V)}$ satisfies the  Davies-Gaffney estimates
\begin{align*}
|\langle e^{t (B^n-V)}f_1, f_2\rangle|\leq\exp\left(-\frac{r^2}{4t}\right)\|f_1\|_{L^2_{c}}\|f_2\|_{L^2_{c}}
\end{align*}
for all $t>0$,  $U_1$,\ $U_2$ open subsets of $(0,+\infty)$, $r:=d(U_1, U_2)=\min\{|x-y|:x\in U_1,y\in U_2\}$ and $f_i$ in $L^2(U_i,y^cdy)$. 
By \cite[Corollary 4.4]{Coulhon-Sikora} and Lemma \ref{Misura palle} we then obtain for  $z\in\Sigma_{\frac\pi 2-\epsilon}$ and $y,\rho>0$
\begin{align*}
|p_V(z,y,\rho)|&\leq C_\epsilon  \frac{1}{\left(Q_c(y,\sqrt{|z|}\right)^{\frac 1 2}\left(Q_c(\rho,\sqrt {|z|}\right)^{\frac 1 2}}\exp\left(-\frac{|y-\rho|^2}{\kappa_\epsilon |z|}\right)\\[1ex]
&\leq C_\epsilon'|z|^{\frac{c+1}2} \left(\frac{y}{\sqrt {|z|}}\right)^{-\frac c 2} \left(1\wedge \frac{y}{\sqrt {|z|}}\right)^{\frac c 2}\left(\frac{\rho}{\sqrt {|z|}}\right)^{-\frac c 2} \left(1\wedge \frac{\rho}{\sqrt {|z|}}\right)^{\frac c 2}\exp\left(-\frac{|y-\rho|^2}{\kappa_\epsilon |z|}\right).
\end{align*}
This is an equivalent form (after modifying the constant in the exponential) of the estimate in the statement, 
by \cite[Lemma 10.2]{MNS-Caffarelli} with $\gamma_1=\gamma_2=-\frac c 2$.\qed

	\begin{os}
		We remark that in  \cite{Coulhon-Sikora}, the authors work in an abstract metric measure space $(M,d, \mu)$ and assume that the heat kernel $p$ associated with a semigroup $e^{-zL}$,  where $L$ is a  non-negative  self-adjoint operator on $L^2(M,d\mu)$, is continuous  with respect to the space variables. In such a case, in fact, 
		\begin{align*}
			\underset{x\in U_1, y\in U_2}{ sup}|p(z,x,y)|
=\mbox{sup}\{\int_{M}e^{-zL}f_1\overline{f_2}\,d\mu,\quad \|f_1\|_{L^1(U_1,d\mu)}=\|f_2\|_{L^1(U_2,d\mu)}=1\}.
		\end{align*}
		In our setting the continuity assumption on $p$ can be avoided since the proofs of \cite[Theorem 4.1, Corollary 4.4]{Coulhon-Sikora}  hold only assuming that  for a.e. $x,y\in M$
		\begin{align*}
			p(z,x,y)=\lim_{s\to 0}\int_{M}e^{-zL}f_1\overline{f_2}\,d\mu
=\lim_{s\to 0} \frac{1}{\mu(B(x,s)\mu(B(y,s)}\int_{B(x,s)\times B(y,s)}p(z,\bar x ,\bar y),d\mu(\bar x)d\mu(\bar y),
		\end{align*}
		where  $f_1=\frac{\chi_{B(x,s)}}{\mu(B(x,s))}$, $f_2=\frac{\chi_{B(y,s)}}{\mu(B(y,s))}$. 
This holds,  outside a set of zero measure, when the measure $\mu$ is doubling, by the Lebesgue differentiation theorem.
	\end{os}

\subsection{A core for $B^n-V$}
We prove that under mild hypotheses the set $$\mathcal {D}= \left\{u \in C_c^\infty ([0,\infty)): u \ {\rm constant \ in \ a \ neighborhood\  of }\  0  \right \}$$ is a core for $B^n-V$ in $L^2_c$. Note that this is true when $V=0$, by Proposition \ref{core}.
We need some elementary lemmas. Unless explicitly stated, we only assume that $0 \leq V\in L^1_{loc}\left(\R^+,y^c\, dy\right)$.

\begin{lem} \label{density}
Assume that $0 \leq V\in L^2_{loc}(\R_+, y^c\, dy).$  Then $D(\mathfrak{a}_V)=H^1_c\cap L^{2}\left(\R^+,Vy^c\, dy \right)$ is dense in $H^1_c$.
\end{lem}
{\sc Proof.} 
By Proposition \ref{core}, $\mathcal D$  is dense in $D(B^n)$ with respect to the graph norm. Moreover, since $V\in L^2_c$ locally, $\mathcal {D}\subset D(\mathfrak{a}_V)$. The claim follows  from the density of $D(B^n)$  in  $H^1_c$.
\qed

\begin{lem}\label{Bounded pert}
Let $u\in H^1_c$ such that $Vu\in L^2_c$. Then $u\in D(B^n)$ if and only if $u\in D(B^n-V)$. Moreover 
$$(B^n-V)u=Bu-Vu.$$
\end{lem}
{\sc Proof.} Let $u\in D(B^n)$. 
Then $u \in D(\mathfrak{a})$ and there exists  $f \in L^2_c$  such that  
$$\mathfrak{a}(u,v)=\int_0^\infty D_y u D_y \overline{v} y^c\, dy=\int_0^\infty f \overline{v}y^c\, dy$$ for every $v\in H^1_c$. 
Setting $g=f+Vu \in L^2_c$ we have
$$\mathfrak{a}_V(u,v)=\int_0^\infty (D_y u D_y \overline{v}+Vu\overline{v}) y^c\, dy=\int_0^\infty (f+Vu) \overline{v}y^c\, dy$$ for every $v\in H^1_c$ and, in particular, for every $v\in D(\mathfrak{a}_V) \subseteq H^1_c$. Therefore $u\in D(B^n-V)$.
Conversely, if $u\in D(B^n-V)$,  
then $u \in D(\mathfrak{a}_V)$ and there exists  $g \in L^2_c$  such that  
$$\mathfrak{a}_V(u,v)=\int_0^\infty (D_y u D_y \overline{v}+Vu\overline{v}) y^c\, dy=\int_0^\infty g\overline{v}y^c\, dy$$ for every $v\in D(\mathfrak{a}_V).$ Setting $f=g-Vu\in L^2_c$ we have that

$$\mathfrak{a}(u,v)=\int_0^\infty f \overline{v}y^c\, dy$$ for every $v\in D(\mathfrak{a}_V)$, hence for every $v\in H^1_c$,  by Lemma \ref{density}.
\qed

\begin{lem} \label{cutoff}
Let $u\in D(B^n-V)$ and $\eta$ be a smooth function such that $\eta=1$ for $0\leq y\leq 1$ and $\eta=0$ for $y\geq 2$. Then $\eta u \in D(B^n-V)$ and 
$$(B^n-V)(\eta u)=\eta (B^n-V) u+ 2 D_y\eta D_y u+ uD_{yy}\eta +cu\frac{D_y \eta}{y}.$$
 \end{lem}
{\sc Proof.} Let $u\in D(B^n-V)$, then $\eta u\in D(\mathfrak{a}_V)$ and,  setting $f=(B^n-V)u$, 
\begin{align*}
\mathfrak{a}_V(\eta u,v)&=\int_0^\infty (D_y (\eta u) D_y \overline{v}+V\eta u\overline{v}) y^c\, dy=\\&
\int_0^\infty (D_yu D_y(\eta \overline{v})+Vu\eta \overline{v}+uD_y\eta D_y\overline{v}-D_y uD_y\eta \overline{v}) y^c\, dy=\\&
-\int_0^\infty \eta f \overline{v} y^c\, dy-\int_0^\infty D_yuD_y\eta \overline{v} y^c\, dy+ \int_0^\infty u D_y\eta  D_y\overline{v} y^c\, dy\\&= -\int_0^\infty \eta f \overline{v} y^c\, dy-\int_0^\infty D_yuD_y\eta \overline{v} y^c\, dy- \int_0^\infty \overline{v} D_y(u D_y \eta   y^c)\, dy\\&=
 -\int_0^\infty \eta f \overline{v} y^c\, dy-2\int_0^\infty D_yuD_y\eta \overline{v} y^c\, dy- \int_0^\infty \overline{v} uD_{yy}\eta   y^c\, dy
- \int_0^\infty\frac{cu}{y} \overline{v} D_{y}\eta   y^c\, dy
\end{align*}
for every $v\in D(\mathfrak{a}_V)$.
\qed

\begin{lem} \label{appr compac}
Let $u\in D(B^n-V)$. Then there exists $(u_k)\subseteq  D(B^n-V)$ with compact support such that $(u_k) \to u$ in $ D(B^n-V)$.
\end{lem}
{\sc Proof.} Let $\eta$ be a smooth function such that $\eta=1$ for $0\leq y\leq 1$ and $\eta=0$ for $y\geq 2$. Setting $\eta_k(y)=\eta\left(\frac{y}{k}\right)$, by Lemma \ref{cutoff}, $u_k=\eta_k u \in D(B^n-V)$ and 
$$(B^n-V)(\eta_k u)=\eta_k (B^n-V) u+ 2 D_y\eta_k D_y u+ uD_{yy}\eta_k +\frac{cu}{y}D_y\eta_k.$$
Then $u_k\to u$, $\eta_k (B^n-V) u\to  (B^n-V) u$ in $L^2_c$ by dominated convergence and, since $D_y \eta_k=0$ in $[0,1]$,  
 $$\left|D_y\eta_k D_y u+u D_{yy}\eta_k +\frac{cu}{y}D_y\eta_k\right|\leq C\left(\frac{|u|}{k}+\frac{|u|}{k^2}+\frac{|D_y u|}{k}\right)\chi_{[k,\infty[}\to 0.$$ 
\qed

Lemma \ref{cutoff} shows that functions with compact support are a core for $B^n-V$. To show that $\mathcal D$ is a core we need more information on the behavior near $y=0$ of functions in the domain of $B^n-V$.

We start by recalling some well-known facts about the modified Bessel functions   $I_{\nu}$ and $K_{\nu}$
which constitute a basis of solutions of the modified Bessel equation
\begin{equation*}
\label{eq.mbessel}
z^2\frac{d^2v}{dz^2}+
z\frac{dv}{dz}-(z^2+\nu^2)v=0,
\quad \Rp z>0.
\end{equation*}
We recall that for $\Rp z>0$ one has 
\begin{equation*} \label{IK}
I_\nu(z)=\left (\frac{z}{2} \right )^\nu \sum_{m=0}^\infty\frac{1}{m!\,\Gamma (\nu+1+m)}\left (\frac{z}{2}\right )^{2m}, \quad K_\nu(z)=\frac{\pi}{2}\frac{I_{-\nu} (z)-I_\nu(z)}{\sin \pi \nu},
\end{equation*}
where limiting values are taken for the definition of $K_\nu$ when $\nu $ is an integer. The basic properties of these functions we need are collected in the following lemma, see e.g., \cite[Sections 9.6 and 9.7]{abramowitz+stegun}.

\begin{lem}\label{behave}
For $\nu> -1$, $I_\nu$ is increasing and  $K_\nu$ is decreasing (when restricted to the positive real half line). Moreover they  satisfy
the following properties if $z \in \Sigma_{\pi/2-\eps}$.
\begin{itemize}
\item[(i)] $I_\nu(z)\neq 0$ for every  $\Rp z>0$. 
\item[(ii)] $I_{\nu}(z)\approx 
\frac{1}{\Gamma(\nu+1)} \left(\frac{z}{2}\right)^{\nu}, \quad \text{as }|z|\to 0,\qquad I_{\nu}(z)\approx \frac{e^z}{\sqrt{2\pi z}}(1+O(|z|^{-1}),\quad \text{as }|z|\to \infty$.\medskip
\item[(iii)] If $\nu\neq 0$,\quad  $K_{\nu}(z)\approx 
\frac{\nu}{|\nu|}\frac{1}{2}\Gamma(|\nu|) \left(\frac{z}{2}\right)^{-|\nu|}, \qquad K_{0}(z)\approx  -\log z, \qquad \text{as }|z|\to 0$
\\[2ex]
$ K_{\nu}(z)\approx \sqrt{\frac\pi{2z}}e^{-z},\quad \text{as }|z|\to \infty$.\medskip
\item[(iv)] $ I_{\nu}'(z)=I_{\nu+1}(z)+\frac{\nu}{z}I_{\nu}(z)$,\quad  $ K_{\nu}'(z)=-K_{\nu+1}(z)+\frac{\nu}{z}K_{\nu}(z)$, for every  $\Rp z>0$.\\[1ex]
\medskip

\end{itemize}

\end{lem}

Note that 
\begin{equation}\label{Asymptotic I_nu}
|I_\nu(z)|\simeq C_{\nu,\epsilon} (1\wedge |z|)^{\nu+\frac 1 2}\frac{e^{Re z}}{\sqrt {|z|}},\qquad z\in \Sigma_{\frac \pi 2-\epsilon}
\end{equation}
for suitable constants $C_{\nu,\epsilon}>0$ which may be different in lower an in the upper estimate. 

The following estimates of  the resolvent operator of $B^n-V$ is a consequence of the domination property stated in Proposition \ref{kernel Vreal}.

\begin{prop} \label{risolventeB^n -V}
Let $c+1>0$ and $\lambda>0$. Then, for every $f\in L^2_c$, 
$$(\lambda -B^n+V)^{-1}f=\int_0^\infty G(\lambda,y,\rho)f(\rho) \rho^c d\rho$$ with
\begin{equation*}  
0\leq G(\lambda,y,\rho)\leq G^n(\lambda,y,\rho)
\end{equation*}
 where \begin{equation}  \label{resolvent}
G^n(\lambda,y,\rho):= \begin{cases}
 y^{\frac{1-c}{2}} \rho^{\frac{1-c}{2}}\,I_{\frac{c-1}{2}}(\sqrt{\lambda}\,y)K_{{\frac{|1-c|}{2}}}(\sqrt{\lambda}\,\rho)\quad y\leq \rho\\[1.5ex]
 y^{\frac{1-c}{2}} \rho^{\frac{1-c}{2}}\,I_{\frac{c-1}{2}}(\sqrt{\lambda}\,\rho)K_{\frac{|1-c|}{2}}(\sqrt{\lambda}\,y)\quad y\geq \rho,
\end{cases}
\end{equation}
is the integral kernel (taken with respect to the measure $\rho^c d\rho$) of the operator $(\lambda-B^n)^{-1}$.
\end{prop}
{\sc Proof.} Writing $(\lambda-B^n+V)^{-1}=\int_0^\infty e^{-\lambda t}e^{t(B^n-V)}dt$ and using property (i) of Proposition \ref{kernel Vreal} we get that
\begin{align*}
 |(\lambda-B^n+V)^{-1}f|\leq (\lambda-B^n)^{-1}|f|,\qquad \lambda>0,\quad  f\in L^2_c.
 \end{align*}
This yields the domination $G(\lambda, y, \rho) \leq G^n(\lambda, y, \rho)$ (the existence of the kernel  follows by \cite[Proposition 1.9]{ArendtBukhalov} as in Proposition \ref{kernel Vreal}). Formula \eqref{resolvent} is proved in \cite[Proposition 2.4]{MNS-Caffarelli}.\qed

We now prove local pointwise estimates for functions in the domain  of $B^n-V$.

\begin{prop} \label{L inf B^n -V}
Let $c+1>0$. Then there exists $C>0$, independent of $V$, such that for every $u\in D(B^n-V)$ and $0<y<1$
\begin{itemize}
\item[(i)] if $-1<c<3$ 
\begin{align*}
|u(y)|\leq C\left(\|u\|_{L^2_c}+\|(B-V)u\|_{L^2_c}\right), 
\end{align*}
\item[(ii)] if $c=3$  
\begin{align*}
|u(y)|\leq C\left(\|u\|_{L^2_c}+\|(B-V)u\|_{L^2_c}\right)|\log y|^{\frac{1}2},
\end{align*}
\item[(iii)] if $c>3$ 
\begin{align*}
|u(y)|\leq C\left(\|u\|_{L^2_c}+\|(B-V)u\|_{L^2_c}\right)y^{\frac{3-c}2}.
\end{align*}
\end{itemize}
\end{prop}

{\sc Proof.} Let $u\in D(B^n-V)$ and $f=u-(B^n-V)u\in L^2_c$ so that $u=(I-B^n+V)^{-1}f$. Let us distinguish between the  following cases and always take $0 <y<1$.
\begin{itemize}
 \item[(i)] If $-1<c<1$,   Lemma \ref{behave} implies that for $y\leq 1$
 \begin{align*}
 G(1,y,\rho)\simeq \begin{cases}
 1,&\quad \rho<1,\\[1ex]
\rho^{-\frac c 2}e^{-\rho},&\quad 1<\rho.
\end{cases}
 \end{align*}
Then  one has
\begin{align*}
|u(y)|&\leq \int_{0}^{\infty}G(1,y,\rho)|f(\rho)|\rho^c d\rho\le C\left (\int_0^1|f(\rho)|\rho^{c}d\rho+ \int_1^\infty  \rho^{-\frac{c}{2}} e^{-\rho}|f(\rho)|\,\rho^{c }d\rho \right )\\[1ex]
 &\leq C\left(\|f\|_{L^2_c((0,1))}+\|\rho^{-\frac c 2}e^{-\rho}\|_{L^2_c((1,\infty))}\|f\|_{L^2_c((1,\infty))}\right)
\leq C\|f\|_{L^2_c}.
\end{align*}
 \item[(ii)] If  $c=1$,    Lemma \ref{behave} gives for $y\leq 1$
 \begin{align*}
 G(1,y,\rho)\simeq \begin{cases}
 |\log y|\leq |\log \rho|,&\quad \rho<y<1,\\[1ex]
|\log \rho|,&\quad y<\rho<1,\\[1ex]
\rho^{-\frac 1 2}e^{-\rho},&\quad 1<\rho.
\end{cases}
 \end{align*}
Then analogously
\begin{align*}
 |u(y)|&
 \le C\left (\int_0^1|\log\rho||f(\rho)|\rho d\rho+ \int_1^\infty  \rho^\frac{1}{2} e^{-\rho}|f(\rho)|d\rho \right )\\[1ex]
 &\leq C\left(\|\log\rho\|_{L^2_c((0,1))}\|f\|_{L^2_c((0,1))}+\|\rho^{-\frac 1 2 }e^{-\rho}\|_{L^2_c((1,\infty))}\|f\|_{L^2_c((1,\infty))}\right)
\leq C\|f\|_{L^2_c}.
\end{align*}
\item[(iii)] Let now   $1<c$. Then   Lemma \ref{behave} implies that for $y\leq 1$
 \begin{align*}
 G(1,y,\rho)\simeq \begin{cases}
 y^{1-c}\leq \rho^{1-c},&\quad \rho<y<1,\\[1ex]
\rho^{1-c},&\quad y<\rho<1,\\[1ex]
\rho^{-\frac c 2}e^{-\rho},&\quad 1<\rho.
\end{cases}
 \end{align*}
If $c<3$ one has 
\begin{align*}
 |u(y)|&
 \le C\left (\int_0^1\rho^{1-c}|f(\rho)|\rho^c d\rho+ \int_1^\infty  \rho^{-\frac{c}{2}} e^{-\rho}|f(\rho)|\rho^cd\rho \right )\\[1ex]
 &\leq C\left(\|\rho^{1-c}\|_{L^2_c((0,1))}\|f\|_{L^2_c((0,1))}+\|\rho^{-\frac c 2 }e^{-\rho}\|_{L^2_c((1,\infty))}\|f\|_{L^2_c((1,\infty))}\right)
\leq C\|f\|_{L^2_c}.
\end{align*}

If $c=3$ then we get
\begin{align*}
 |u(y)|&
 \le C\left (y^{-2}\int_0^y|f(\rho)|\rho^3 d\rho+\int_y^1\rho^{-2}|f(\rho)|\rho^3 d\rho+ \int_1^\infty  \rho^{-\frac{3}{2}} e^{-\rho}|f(\rho)|\rho^3d\rho \right )\\[1ex]
 &\leq C\|f\|_{L^2_c}\left(y^{-2}\left(\int_0^y\rho^3 d\rho\right)^{\frac 1 2 }+\left(\int_y^1\rho^{-4}\rho^3 d\rho\right)^{\frac 1 2 }+\|\rho^{-\frac 3 2 }e^{-\rho}\|_{L^2_c((1,\infty))}\right)\\[1ex]
 &\leq C\|f\|_{L^2_c}\left(1+|\log y|^{\frac 1 2 }\right)
\end{align*}

and finally if $c>3$ 
\begin{align*}
 |u(y)|&
 \le C\left (y^{1-c}\int_0^y|f(\rho)|\rho^c d\rho+\int_y^1\rho^{1-c}|f(\rho)|\rho^c d\rho+ \int_1^\infty  \rho^{-\frac{c}{2}} e^{-\rho}|f(\rho)|\rho^cd\rho \right )\\[1ex]
 &\leq C\|f\|_{L^2_c}\left(y^{1-c}\left(\int_0^y\rho^c d\rho\right)^{\frac 1 2 }+\left(\int_y^1\rho^{2-2c}\rho^c d\rho\right)^{\frac 1 2 }+\|\rho^{-\frac c 2 }e^{-\rho}\|_{L^2_c((1,\infty))}\right)\\[1ex]
 &\leq C\|f\|_{L^2_c}y^{\frac{3-c} 2}.
\end{align*}
\end{itemize} 
\qed

We can now show that,  under stronger assumptions,  the potential term $V$  can be seen as a  perturbation of $B^n$ near 0, that is $Vu \in L^2_c$ for every $u \in D(B_n)$ having compact support. In particular we prove that $\mathcal D$ is  a core for $B^n-V$.
\begin{prop} 
\label{core B-V}
Let $c+1>0$ and assume that
\begin{itemize}
\item[(i)]
$c<3$ and $V\in L^2_{loc}\left(\R^+,y^cdy\right)$ or
\item[(ii)]
$c=3$ and $V|\log y|^{\frac1 2 }\in L^2_{loc}\left(\R^+,y^cdy\right)$ or
\item[(iii)]$c>3$ and $Vy^{\frac{3-c}2}\in L^2_{loc}\left(\R^+,y^cdy\right)$.
\end{itemize}
 If $\mathcal{C}_r:=\{u\in L^2_c: \mbox{supp }u\subseteq [0,r]\}$ then 
$
 D(B^n-V)\cap \mathcal{C}_r=D(B^n)\cap \mathcal{C}_r$
 with equivalence of norms
 \begin{align*}
 \|u\|_{D(B^n-V)}&\simeq \|u\|_{D(B^n)},\quad \forall u\in D(B^n)\cap \mathcal{C}_r.
 \end{align*}
Finally,  $$\mathcal {D}= \left\{u \in C_c^\infty ([0,\infty)): u \ {\rm constant \ in \ a \ neighborhood\  of }\  0  \right \}$$
is a core for $B^n-V$. 
\end{prop}
{\sc{Proof.}} Let $u\in \mathcal{C}_r$. Then the hypotheses on $V$ and Proposition \ref{L inf B^n -V} imply that $Vu\in L^2_c$ and $\|Vu\|_{L^2_c}\leq C\|u-(B-V)u\|_{L^2_c}$. Then  by  Lemma \ref{Bounded pert} $u\in D(B^n-V)$ if and only if $u\in D(B^n)$. This show the equality $ D(B^n-V)\cap \mathcal{C}_r=D(B^n)\cap \mathcal{C}_r$. Using Proposition \ref{L inf B^n -V} again  we also  have $\|Vu\|_{L^2_c}\leq C_1\|u-Bu\|_{L^2_c}$  for any $u\in D(B^n)\cap \mathcal{C}_r$, which proves the  equivalence of the graph norms. 
Finally,  let  $u\in D(B^n-V)$. We have to prove that $u$ can be approximated in the graph norm with functions belonging to $\mathcal D$. Using Lemma \ref{appr compac} we may suppose, without any loss of generality, that $\mbox{supp\,}u\subseteq (0,r)$. Then, by Proposition \ref{core},  there exist $(u_n) \subset \mathcal D$ such that $u_n\to u$ in the graph  norm $\|\cdot\|_{D(B^n)}$. We may also assume, after multiplying by a suitable cut-off function, that $\mbox{supp\,}u_n\subseteq (0,2r)$ for every $n$.
Then the previous point implies that $Vu_n\to Vu$ in $L^2_c$, too.\qed

\section{The operator $y^\alpha B^n-V$ in $L^2_{c-\alpha}$}

We consider now  for $c \in\R$, $\alpha<2$,  and $0\leq V\in L^1_{loc}\left(\R^+,y^{c-\alpha}\, dy\right)$ the  operator 
\begin{equation*} 
y^\alpha B^n-V =y^\alpha\left(D_{yy}+\frac{c}{y}D_y  \right)-V
\end{equation*}
in the space $L^2_{c-\alpha}$. As in Section \ref{sec y^aB}   we use the isometry  $T_{-\frac \alpha 2}u(y)=\left|1-\frac{\alpha}{2}\right |^{\frac 1p}u(y^{1-\frac{\alpha}{2}})$,
$$T_{-\frac \alpha 2}:L^2_{\tilde c}\to L^2_{c-\alpha},\qquad \tilde c=\frac{c-\frac\alpha 2}{1-\frac\alpha 2},$$
under which   $y^\alpha B-V$ becomes  similar to 
\begin{align*}
T_{-\frac \alpha 2}^{-1} \Big(y^{\alpha}B-V\Big)T_{-\frac \alpha 2}=\left(1-\frac\alpha 2\right)^2\left(\tilde B-\tilde V\right)
\end{align*}
where $\tilde B=D_{yy}+\frac{\tilde c}{y}D_y $ and $\tilde V(y)=(1-\frac \alpha 2)^{-2}V\left(y^{\frac 2 {2-\alpha}}\right)\in L^1_{loc}\left(\R^+,y^{\tilde c}\, dy\right)$.

Defining
\begin{equation*} 
D(y^\alpha B^n-V):=T_{-\frac \alpha 2}\left (D(\tilde{B^n}-\tilde V) \right )
\end{equation*}
one obtains that, when $c>-1+\alpha$, $y^\alpha B^n-V$ generates a contractive analytic semigroup $\left\{e^{z(y^\alpha B^n-V)}:\ z\in\C_+\right\}$  in $L^2_{c-\alpha}$ which satisfies
\begin{align}\label{Op V equiv}
e^{z\left(y^{\alpha}B-V\right)}=T_{-\frac \alpha 2}\left(e^{z\left(1-\frac\alpha 2\right)^2\left(\tilde B-\tilde V\right)}\right)T_{-\frac \alpha 2}^{-1} .
\end{align}

We state the properties  obtained so far, together with a density result which is a restating of Proposition \ref{core B-V} under the isometry $T_{-\frac \alpha 2}$.
\begin{prop} 
\label{core yalpha B-V}
Let $c+1-\alpha>0$ and $0 \leq V\in L^1_{loc}\left(\R^+,y^{c-\alpha}\right)$. Then the operator $y^\alpha B^n-V$ generates a contractive analytic semigroup   in $L^2_{c-\alpha}$. If, in addition, 
\begin{itemize}
\item[(i)]
$c<3-\alpha$ and $V\in L^2_{loc}\left(\R^+,y^{c-\alpha}\right)$ or
\item[(ii)]
$c=3-\alpha$ and $V|\log y|^{\frac 1 2 }\in L^2_{loc}\left(\R^+,y^{c-\alpha}\right)$ or 
\item[(iii)]$c>3-\alpha$ and $Vy^{\frac{3-c-\alpha}2}\in L^2_{loc}\left(\R^+,y^{c-\alpha}\right)$
\end{itemize}
 then 
 $$\mathcal {D}= \left\{u \in C_c^\infty ([0,\infty)): u \ {\rm constant \ in \ a \ neighborhood\  of }\  0  \right \}$$
is a core for $y^\alpha B^n-V$ in $L^2_{c-\alpha}$. 
\end{prop}

\begin{os}\label{os core yalpha B-V}
If $V(y)=y^\alpha$, then $V$ always satisfies (ii) and (iii) when $c \geq 3-\alpha$. Instead, if $c<3-\alpha$, we need $c+1-|\alpha|>0$.

\end{os}

Let $\mathfrak{a}_{\tilde V}$ be the form in $L^2_{\tilde c}$, defined in \eqref{form v},  associated to  $\tilde B^n-\tilde V$. In $L^2_{c-\alpha}$ we introduce  the form $\mathfrak{a}_{\alpha, V} $ which is the image  of $\mathfrak{a}_{\tilde V}$ under the isometry $T_{0,-\frac \alpha 2}$, that is 
\begin{align}\label{form alpha V}
\nonumber	&\mathfrak{a}_{\alpha, V}(u,v):=\mathfrak{a}_{ \tilde V}\left(T_{-\frac \alpha 2}^{-1}u,T_{-\frac \alpha 2}^{-1}v\right)=
	\int_{\R_+}
	\left(y^\alpha D_y u D_y \overline{v}+Vu\overline{v}\right)y^{c-\alpha}\,dy,\\[1ex]
	&D(\mathfrak{a}_{\alpha, V}):=T_{-\frac \alpha 2}D(\mathfrak{a}_{\tilde V})=\left\{u\in L^2_{c-\alpha}:u'\in L^2_{c}\right\}\cap L^{2}\left(\R^+,Vy^{c-\alpha}\, dy\right).
\end{align}
To keep consistency of notation we often write $\mathfrak{a}_{0,V}=\mathfrak{a}_{V}$.  By construction, $y^\alpha B^n-V$  is the operator  associated to the form $\mathfrak{ a}_{\alpha,V}$ in $L^2_{c-\alpha}$
\begin{align*}
	 D(y^\alpha B^n-V)&=\{u \in D(\mathfrak{a}_{\alpha,V}): \exists  f \in L^2_{c-\alpha} \quad {\rm such\ that}\ \\[1ex]
	&\hspace{25ex}   \mathfrak{a}_{\alpha,V}(u,v)=\int_0^\infty f \overline{v}y^{c-\alpha}\, dy
 \	{\rm for\ every}\ v\in D(\mathfrak{a}_{\alpha,V})\},\\
	y^\alpha B^nu-Vu&=-f.
\end{align*}
The next lemma, which follows from the considerations above,  will be used later to relate the resolvents  of $y^\alpha B^n- y^\alpha$  and  $B^n- y^{-\alpha}$.

\begin{lem}\label{os form coinc}
	Let $\mathfrak{a}_{\alpha, y^\alpha}$ and  $\mathfrak{a}_{y^{-\alpha}}$ be the sesquilinear forms associated  respectively to the operator $y^\alpha B^n- y^\alpha$   in $L^2_{c-\alpha}$ and  $B^n- y^{-\alpha}$ in  $L^2_c$. Then 
	\begin{align*}
	\mathfrak{a}_{\alpha, y^\alpha }(u,v)=
	\int_{\R_+}
	\left( D_y u D_y \overline{v}+u\overline{v}\right)y^{c}\,dy,\qquad\mathfrak{a}_{ y^{-\alpha} }(u,v)=
	\int_{\R_+}
	\left( D_y u D_y \overline{v}+y^{-\alpha}u\overline{v}\right)y^{c}\,dy.
\end{align*}
on the common form domain
\begin{align*}
	D(\mathfrak{a}_{\alpha, y^\alpha})=D(\mathfrak{a}_{ y^{-\alpha}})=\left\{u\in L^2_{c-\alpha}\cap L^2_c:u'\in L^2_{c}\right\}
\end{align*}
\end{lem}

Note that the above operators act in different Hilbert spaces; in particular their domains are different. However, the form domains coincide.

\section{The operator $y^\alpha B^n-V$ in $L^p_{m}$}

Here we investigate properties of  $y^{\alpha}B-V$, $\alpha<2$,  in $L^p_m$ when
$0<\frac{m+1}p<c+1-\alpha$. 

We introduce  the  family of integral operators $\left(S_{\alpha}^\beta(t)\right)_{t>0}$  on $L^p_m$
$$
S_{\alpha}^{\beta}(t)f(y):=t^{-\frac {1} 2}\,\int_{\R_+}  \left (\frac{\rho}{t^{\frac{1}{2-\alpha}}}\wedge 1 \right)^{-\beta+\frac\alpha 2}
\exp\left(-\frac{|y^{1-\frac\alpha 2}-\rho^{1-\frac\alpha 2}|^2}{\kappa t}\right)f(\rho) \rho^{-\frac\alpha 2}\,d\rho
$$
and note that 
$$
S_{\alpha}^{\beta}(t)=T_{-\frac \alpha 2}\circ S_0^{\tilde\beta}(t)\circ T_{-\frac \alpha 2}^{-1},\qquad \tilde \beta=\frac{\beta-\frac\alpha 2}{1-\frac\alpha 2}.
$$
As usual  $T_{-\frac \alpha 2}u(y)=\left|1-\frac{\alpha}{2}\right |^{\frac 1p}u(y^{1-\frac{\alpha}{2}})$ is an isometry from $L^p_{\tilde m}$ onto $L^p_m$, $\tilde m=\frac{m+\frac\alpha 2}{1-\frac\alpha 2}$.
Here $\kappa$ is a positive constant but we omit the dependence on it.  The following result has been proved for $\alpha=0$ in \cite{MNS-Caffarelli}.   
 \begin{lem} \label{gen R-bound sempre}
 Let   $m\in\R$ and let  $p\in(1,\infty)$ such that $0<\frac{m+1}p<1-\alpha-\beta$. The families  $\left(S_\alpha^{\beta}(t)\right)_{t\geq 0}$  and $\{\Gamma(\lambda)=\int_0^\infty \lambda e^{-\lambda t}S_\alpha^{\beta}(t)\, dt, \quad \lambda>0\}
$ are  $\mathcal{R}$-bounded in $L^p_m$.
\end{lem}
{\sc{Proof.}} Since the  $\mathcal{R}$-boundedness  is preserved under isometries, from $S_{\alpha}^{\beta}(t)=T_{-\frac \alpha 2}\circ S_0^{\tilde \beta}(t)\circ T_{-\frac \alpha 2}^{-1}$  we may assume that $\alpha=0$ (note that $0<\frac{m+1}p<-\beta+1-\alpha$ is equivalent to $0<\frac{\tilde m+1}p<-\tilde \beta+1$). The first result is then a consequence of \cite[Theorem 7.7]{MNS-Caffarelli}. The family
$$
\Gamma (\lambda)=\int_0^\infty \lambda e^{-\lambda t}S_\alpha^{\beta}(t)\, dt, \quad \lambda>0
$$
is $\mathcal R$-bounded by \cite[Corollary 2.14]{KunstWeis}.\qed

We can now prove our main results for the operator  $y^{\alpha}B-V$. 

\begin{teo}\label{Gen pert R}
Let  $0\leq V\in L^1_{loc}\left(\R^+,y^{c-\alpha}\, dy\right)$. For any $p\in(1,\infty)$ such that $0<\frac{m+1}p<c+1-\alpha$, the semigroup $e^{z(y^\alpha B^n-V)}$ initially defined   on $L^2_{c-\alpha}$ extends to  a bounded analytic semigroup on $L^p_m$ of angle $\pi/2$ which consists of integral operators. Moreover the generated semigroup  has maximal regularity and the following properties hold.
\begin{itemize}
  \item[(i)] For every $\epsilon>0$ there exist $C=C(\epsilon,\alpha)>0$ (independent of $V$) such that
     \begin{align*}
   \left|e^{z(y^\alpha B^n-V)}f\right|\leq CS_\alpha^{-c}(|z|)|f|,\quad f\in L^p_{m},\quad |\arg z|<\frac\pi 2-\epsilon.
   \end{align*}
\item[(ii)] For every $\epsilon >0$  the families of operators  
\begin{align*}
\left\{e^{z(y^\alpha B^n-V)}:\; z\in \Sigma_{\frac \pi 2-\epsilon},\; 0\leq V\in L^1_{loc}\left(\R^+,y^{c-\alpha}\right)\right\},\\[1ex]
\left\{\lambda\left(\lambda-y^\alpha B^n+V\right)^{-1}:\; \lambda \in \Sigma_{\pi-\epsilon}:\; 0\leq V\in L^1_{loc}\left(\R^+,y^{c-\alpha}\right)\right\}
\end{align*}
 are $\mathcal R$-bounded in $L^p_m$.
  \end{itemize}  
\end{teo}  
{\sc Proof. } 
By  Proposition \ref{kernel V} and \eqref{Op V equiv}, (i) holds for any $f\in L^2_{c-\alpha}$. The boundedness of $e^{z(y^\alpha B^n-V)}$ in $L^p_m$ follows from the previous lemma and (i) extends to $L^p_m$. The semigroup law is inherited  from  $L^2_{c-\alpha}$ via a density argument and we have only to prove the  strong continuity at $0$. Using the isometry $T_{-\frac \alpha 2}$, we may suppose that $\alpha=0$. Let $f, g \in C_c^\infty (0, \infty)$. Then as $z \to 0$, $z \in \Sigma_{\frac\pi 2-\epsilon}$, 
$$
\int_0^\infty (e^{z(B^n-V)}f)\, g\, y^m dy=\int_0^\infty (e^{z(B^n-V)}f) \,g\, y^{m-c}y^c  dy \to \int_0^\infty fgy^{m-c}y^cdy  =\int_0^\infty fgy^{m}dy,
$$
by the strong continuity of $e^{z(B^n-V)}$ in $L^2_c$. By density and uniform boundedness of the family $(e^{z(B^n-V)})_{ z\in \Sigma_{\frac\pi 2-\epsilon}}$ this holds for every $f \in L^p_m$, $g \in L^{p'}_m$. The semigroup is then weakly continuous, hence strongly continuous.  

  The $\mathcal R$-boundedness of e$^{z(y^\alpha B^n-V)}$ follows then  by domination from Lemma \ref{gen R-bound sempre}, see Corollary \ref{domination}. To prove the $\mathcal R$-boundedness of the resolvent family, for $\lambda \in \Sigma_{\pi-\epsilon}\setminus\{0\}$ let  $\theta=\frac{|\mbox{arg}\lambda|}{\mbox{arg}\lambda}\left(\frac\pi 2-\frac\epsilon 2\right)$ so that $\mu:=e^{-i\theta}\lambda \in\Sigma_{\frac\pi 2-\frac\epsilon 2}$. Then
\begin{align*}
   \left|\lambda\left(\lambda-y^\alpha B^n+V\right)^{-1}f\right|&=  \left|\mu\left(\mu-e^{-i\theta}(y^\alpha B^n-V)\right)^{-1}f\right|= \left|\int_0^\infty \mu e^{-\mu t}e^{-i\theta t(y^\alpha B^n-V)}f\,dt\right|\\[1ex]
&\leq  C\int_0^\infty |\mu| e^{-\Rp\mu t}S^{-c}_\alpha(t)|f|\,dt\leq C\int_0^\infty |\lambda| e^{-|\lambda|\sin\frac\epsilon 2 t}S^{-c}_\alpha(t)|f|\,dt . 
 \end{align*}
The $\mathcal R$-boundedness of the second family in (ii) now follows from \cite[Corollary 2.14]{KunstWeis} and the maximal regularity of the semigroup from  Theorem \ref{MR}.\qed

In our investigation of  degenerate Nd problems, see \cite{MNS-CompleteDegenerate}, we need  also a weaker version of the result above for potentials having non-negative real part. We formulate it in the next proposition.

\begin{prop}\label{Gen pert C}
Let $V\in L^1_{loc}\left(\R^+,y^{c-\alpha}\, dy\right)$ be a   potential having non-negative real part. Then for any $1<p<\infty$ such that $0<\frac{m+1}p<c+1-\alpha$,  $y^\alpha B^n-V$ generates a $C_0$-semigroup on $L^p_m$. The generated semigroup consists of integral operators and the following estimates hold 
     \begin{align*}
   \left|e^{t(y^\alpha B^n-V)}f\right|&\leq e^{ty^\alpha B^n}|f|,\hspace{13ex} f\in L^p_m,\quad t\geq 0
   \end{align*}
In particular the families of operators
\begin{align*}
	\left\{e^{t(y^\alpha B^n-V)}:\; t\geq 0,\;  V\in L^1_{loc}\left(\R^+,y^{c-\alpha}\right), \ \Rp V\geq 0\right\},\\[1ex]
	\left\{\lambda\left(\lambda-y^\alpha B^n+V\right)^{-1}:\; \lambda >0,\; V\in L^1_{loc}\left(\R^+,y^{c-\alpha}\right), \ \Rp V\geq 0\right\}
\end{align*}
are $\mathcal R$-bounded in $L^p_m$.
\end{prop}
{\sc{Proof.}}
Using the isometry $T_{0,-\frac \alpha 2}$ we may assume that $\alpha=0$. Let us treat first the  symmetric case  in $L^2_c$.  The generation results can be proved as  in  Proposition \ref{kernel Vreal} (where we assumed $V\geq 0$). If $\mathfrak a$ is the  form associated with $ B^n$, then  $ B^n-V$ is associated to ${\mathfrak a}_V:=\mathfrak a(u,v)+\langle Vu,v\rangle_{L^2_c}$ and,  by the standard theory on sesquilinear forms,  $ B^n-V$  generates a $C_0$-semigroup on $L^2_c$. 

The domination properties follow from  \cite[Theorem 2.21]{Ouhabaz}. Let $u,v\in D\left({\mathfrak a}_V\right)=D(\mathfrak{a})\cap L^{2}\left(\R^+,|V|y^c\, dy\right)$ such that $u\bar v\geq 0$. Since $e^{t B^n}$ is positive one has $\Rp \mathfrak a(u,v)\geq \mathfrak a(|u|,|v|)$. Moreover 
\begin{align*}
\Rp {\mathfrak a}_V(u,v)=\Rp \mathfrak a(u,v)+\int_0^\infty  \Rp V\, u\bar v \,y^cdy\geq \Rp\mathfrak a(|u|,|v|)
\end{align*}
which by \cite[Theorem 2.21]{Ouhabaz} again  implies the stated domination of the generated semigroups (one easily verifies that $D\left({\mathfrak a}_V\right)$ is an ideal of $ D(\mathfrak{a})$ since this last is an ideal in itself, by the positivity of $e^{tB^n}$, see \cite[Proposition 2.20]{Ouhabaz}). The extrapolation on $L^p_m$ follows as in Theorem \ref{Gen pert R}.  The domination of the resolvent is a straightforward consequence of that of the semigroup. The $\mathcal R$-boundedness of the semigroup   follows  by domination from the  $\mathcal R$-boundedness of $(e^{t B^n})_{ t \geq 0}$ proved in Theorem \ref{Gen pert R}.  The $\mathcal R$-boundedness of the resolvent follows as in Theorem \ref{Gen pert R}. \qed

\section{The operator $y^\alpha B^n-y^\alpha$ }\label{case V yalpha}
 We end the paper by thoroughly  investigating  the special case $V(y)=y^\alpha$, keeping $\alpha<2$. We prove, in particular, that the domain of $y^\alpha B-V$ is $D(y^\alpha B)\cap D(V)$, under slightly more restrictive hypotheses than those of Theorem \ref{Gen pert R}.

As explained in the Introduction this case  plays   a crucial role in \cite{MNS-CompleteDegenerate} in  the investigation  of the  degenerate operators
	\begin{align*}
		\mathcal L =y^{\alpha_1}\Delta_{x} +y^{\alpha_2}\left(D_{yy}+\frac{c}{y}D_y  -\frac{b}{y^2}\right),\qquad\alpha_1, \alpha_2 \in\R
	\end{align*}
	in the spaces $L^p\left(\R^{N+1}_+, y^mdxdy\right)$. 
 In particular 	we prove in Propositions \ref{rb y^aResolv} and  \ref{Lem Mult Resolv} that  the multipliers 
 \begin{align*}
 	\xi \in \R^N \to N_{\lambda}(\xi)&=\lambda (\lambda-y^\alpha By+y^\alpha|\xi|^2)^{-1},\\[1ex]
 	\xi \in \R^N \to M_{\lambda}(\xi)&=y^\alpha|\xi|^2(\lambda-y^\alpha By+y^\alpha|\xi|^2)^{-1}
 \end{align*}
 satisfy the hypothesis of  Theorem \ref{mikhlin}.

We start with the following Lemma.
\begin{lem}\label{Lemma Rb 3}
	Assume that   $c+1>0$ and  $c+1-\alpha>0$, that is  $B^n$ generates a $C_0$-semigroup in $L^2_{c}$ and  $y^\alpha B^n$ generates a $C_0$-semigroup in $L^2_{c-\alpha}$. If  $ \lambda \in \mathbb C^+$ and $\mu> 0$,  then 
	\begin{align*}
	\left(\lambda-y^\alpha B^n+\mu y^\alpha\right)^{-1}f=\left(\mu- B^n+\frac{\lambda}{y^\alpha}\right)^{-1}\left(\frac{f}{y^\alpha}\right),\qquad \forall f\in C^\infty_c((0,\infty)).
\end{align*}	
\end{lem} 
{\sc{Proof.}} Under the assumptions  $y^\alpha B^n-\mu y^\alpha$ and $B^n-\lambda y^{-\alpha}$  generate a semigroup  on $L^2_{c-\alpha}$ and $L^2_c$, respectively, see Theorem \ref{Gen pert R}. Since $\Rp \lambda>0$, $\mu>0$ both resolvents are well defined but map to different spaces.

Let  $\mathfrak{a}_{\alpha,\mu y^\alpha}$ , $\mathfrak{a}_{\lambda y^{-\alpha}}$ be the forms  associated   to $y^\alpha B^n-\mu y^\alpha$   in $L^2_{c-\alpha}$ and  $B^n-\lambda y^{-\alpha}$ in  $L^2_c$
\begin{align*}
	\mathfrak{a}_{\alpha,\mu y^\alpha }(u,v)=
	\int_{\R_+}
	\left( D_y u D_y \overline{v}+\mu u\overline{v}\right)y^{c}\,dy,\qquad\mathfrak{a}_{ \lambda y^{-\alpha} }(u,v)=
	\int_{\R_+}
	\left( D_y u D_y \overline{v}+\lambda y^{-\alpha}u\overline{v}\right)y^{c}\,dy.
\end{align*}
By Lemma \ref{os form coinc} they are defined on the common domain
\begin{align*}
	\mathcal{F}:=\left\{u\in L^2_{c-\alpha}\cap L^2_c:u'\in L^2_{c}\right\}
\end{align*}

 Given   $f\in C^\infty_c((0,\infty))$ let  $u:= \left(\mu- B^n+\frac{\lambda}{y^\alpha}\right)^{-1}\left(\frac{f}{y^\alpha}\right)$. In order to prove  that the equality  $u=\left(\lambda-y^\alpha B^n+\mu y^\alpha\right)^{-1}f$ holds,  we have to show that $u\in \mathcal{F}$ and that for every $v\in \mathcal F$, $u$ satisfies the weak equality 
 \begin{align}\label{Rb 3 eq1}
\int_0^\infty f \overline{v}y^{c-\alpha}\, dy&=\int_0^\infty \lambda u \overline{v}y^{c-\alpha}\, dy+\mathfrak{a}_{\alpha,\mu y^\alpha}(u,v)=\int_0^\infty (\lambda y^{-\alpha}u \overline{v}+ D_y u D_y \overline{v}+\mu u\overline{v})y^{c}\, dy.
 \end{align}
 By construction  $u$ is in the domain  of $B^n-\lambda y^{-\alpha}$ which is contained in $\mathcal F$ and satisfies  
 \begin{align*}
	\int_0^\infty \frac{f}{y^\alpha} \overline{v}y^{c}\, dy&=\int_0^\infty \mu u \overline{v}y^{c}\, dy+\mathfrak{a}_{\alpha,\lambda y^{-\alpha}}(u,v)=\int_0^\infty (\mu u \overline{v}+ D_y u D_y \overline{v}+\lambda y^{-\alpha} u\overline{v})y^{c}\, dy,
\end{align*}
which is the same as  \eqref{Rb 3 eq1}. \qed

%
%
%

In the next results we relate the resolvent of $y^\alpha B^n- y^\alpha$ with that of $B^n-\frac{1}{y^\alpha}$.  We shall assume both the conditions  $0<\frac{m+1}p<c+1-\alpha$ and $-\alpha<\frac{m+1}p<c+1-\alpha$ (that is $\alpha^-<\frac{m+1}p<c+1-\alpha$). The first guarantees that $y^\alpha B^n$ is a generator in $L^p_{m}$ and the second that $B^n$ is a generator in $L^p_{m+\alpha p}$.

\begin{cor}\label{Cor Rb}
	Assume that  $\alpha^-<\frac{m+1}p<c+1-\alpha$. If $ \lambda \in \mathbb C^+$ and $\mu> 0$, then 
	\begin{itemize}
		\item[(i)] for every $f\in L^p_m$ 
		\begin{align*}
			\left(\lambda-y^\alpha B^n+\mu y^\alpha\right)^{-1}f=\left(\mu- B^n+\frac{\lambda}{y^\alpha}\right)^{-1}\left(\frac{f}{y^\alpha}\right)\in L^p_{m+\alpha p}\cap L^p_m;
		\end{align*}
		
		\item[(ii)] the operator 	$ y^\alpha\left(\lambda-y^\alpha B^n+\mu y^\alpha\right)^{-1}$ is  bounded in $L^p_m$;
\item[(iii)] the operator 	$ \frac{1}{y^\alpha}\left(\mu- B^n+\frac{\lambda}{y^\alpha}\right)^{-1}$
		is  bounded in $L^p_{m+\alpha p}$.
	\end{itemize}
	
\end{cor} 
{\sc{Proof.}} Equality (i) is proved in Lemma \ref{Lemma Rb 3} for any $f\in C^\infty_c((0,\infty))$. Since  $\left(\lambda-y^\alpha B^n+\mu y^\alpha\right)^{-1}$ is bounded form $L^p_m$ into itself and $\left(\mu- B^n+\frac{\lambda}{y^\alpha}\right)^{-1}\left(\frac{\cdot}{y^\alpha}\right)$ is bounded from $L^p_m$ to $L^p_{m+\alpha p}$,  by density, (i) holds for every $f \in L^p_m$.
Parts (ii), (iii) are consequence of (i).\\\qed
\smallskip
In the next propositions we prove the boundedness of the multipliers $N_{\lambda}$ and $M_{\lambda}$. We start with 
$M_{\lambda}$, used in  \cite{MNS-CompleteDegenerate}  to  characterize the domain of  $\mathcal L=y^\alpha(\Delta_x+B_y)$.

\begin{prop}\label{rb y^aResolv}
Assume that  $\alpha^-<\frac{m+1}p<c+1-\alpha$ and let for $\lambda \in \C^+$, $\xi \neq 0$
$$M_\lambda (\xi)=|\xi|^2 y^\alpha\left(\lambda-y^\alpha B^n+|\xi|^2 y^\alpha\right)^{-1}\in\mathcal {B}(L^p_m).$$
Then the family
$\left \{|\xi|^{|\beta|}D^\beta_\xi( M_{\lambda})(\xi): \xi\in \R^{N}\setminus\{0\}, \ |\beta | \leq N ,\lambda \in \C^+\right \}$
is $\mathcal{R}$-bounded in $L^p_m$.


\end{prop}
{\sc{Proof.}} Let $m_\lambda (\mu)=\mu y^\alpha  \left(\lambda -y^\alpha B^n+\mu y^\alpha \right)^{-1}$, $\mu>0$. 

Using Lemma \ref{mikhlin Comp} it suffices to show that the family$\left \{\mu ^{k}D^k_\mu ( m_{\lambda})(\mu): \mu>0, \ k \leq N ,\lambda \in \C^+\right \}$
is $\mathcal{R}$-bounded in $L^p_m$.

The map $Tf=f/y^\alpha$ is an isometry of $L^p_m$ onto $L^p_{m+\alpha p}$ and by Corollary \ref{Cor Rb}, $$m_\lambda (\mu)=T^{-1}\mu \left(\mu-B^n+\frac{\lambda}{y^\alpha}\right)^{-1}T.$$ 

The family 
$$\left \{\mu ^{k}D^k_\mu( \Gamma_{\lambda})(\mu): \mu>0, \ k \leq N ,\lambda \in \C^+\right \},\qquad \Gamma_\lambda (\mu)=\mu \left(\mu-B^n+\frac{\lambda}{y^\alpha}\right)^{-1}
$$
 is  $\mathcal{R}$-bounded in $L^p_{m+\alpha p}$. Indeed, 
 $$
\Gamma_\lambda (\mu)=\int_0^\infty \mu e^{-\mu t}e^{t(B^n-\frac{\lambda}{y^\alpha})}\, dt
$$
and  $\left\{e^{t(B^n-\frac{\lambda}{y^\alpha})}:\,  t \geq 0, \  \lambda \in \C^+\right \}$ is $\mathcal R$-bounded in $L^p_{m+\alpha p}$,  by Theorem \ref{Gen pert C}.  The  $\mathcal R$-boundedness of the derivatives follows either by the resolvent equation or  by differentiating the last equation under the integral and using \cite[Corollary 2.14]{KunstWeis}. In fact, if  $h(\mu,t)=\mu e^{-\mu t}$, then
$$\mu^k\int_0^\infty |D^k_\mu h(\mu,t)|dt  \le C_k, \quad  \mu>0.
$$ 
\qed

Next we deal with $N_{\lambda}$ which is crucial in \cite{MNS-CompleteDegenerate} for the proof   that $\mathcal L=y^\alpha(\Delta_x+B_y)$ generates an analytic semigroup.
%

\begin{prop}\label{Lem Mult Resolv}
Assume that  $\alpha^-<\frac{m+1}p<c+1-\alpha$ and let for $\lambda \in \C^+$, $\xi \neq 0$
$$N_{\lambda}(\xi)=(\lambda -y^\alpha B^n+|\xi|^2y^\alpha)^{-1}\in \mathcal B\left(L^p_m\right).$$
Then  the family 
$$\left \{|\xi|^{|\beta|}D^\beta_\xi(\lambda  N_{\lambda})(\xi): \xi\in \R^{N}\setminus\{0\}, \ |\beta | \leq N ,\lambda \in \C^+\right \}$$
is $\mathcal{R}$-bounded in $L^p_m$.
\end{prop}
{\sc{Proof.}}
 For $\mu>0$ let $n_\lambda (\mu)= \left(\lambda -y^\alpha B^n+\mu y^\alpha \right)^{-1}$.  Using Lemma \ref{mikhlin Comp} we have to  show that  
the family 
\begin{align}\label{Lem Mult Resolv eq 5}
	\left \{\mu ^{k}D^k_\mu ( n_{\lambda})(\mu): \mu>0, \ k \leq N ,\ \lambda \in \C^+\right \}
\end{align}
is $\mathcal{R}$-bounded in $L^p_m$.

Theorem \ref{Gen pert R} with $V(y)=\mu y^\alpha$ and Proposition \ref{rb y^aResolv} imply  that the families 
\begin{align}\label{Lem Mult Resolv eq 3}
	\left \{\lambda n_{\lambda}(\mu): \mu>0,\ \lambda \in \C^+\right \},\qquad \left \{\mu y^\alpha  n_{\lambda}(\mu): \mu>0, \ \lambda \in \C^+\right \}
\end{align}
are  $\mathcal R$-bounded in $L^p_m$. 

We have that  $n_{\lambda}(\cdot)\in C^1\left(\R_+,\mathcal B\left(L^p_m\right)\right)$ and
\begin{align}\label{Lem Mult Resolv eq 1}
D_\mu (n_{\lambda}(\mu))=- n_{\lambda}(\mu) y^\alpha n_{\lambda}(\mu).
\end{align}
Indeed setting $A=\lambda -y^\alpha B^n_y$, $V=y^\alpha$ we have 
\begin{align*}
\frac{n_{\lambda}(\mu+h)-n_{\lambda}(\mu)}{h}&= \frac{(A+(\mu+h)V)^{-1}-(A+\mu V)^{-1}}{h}\\[1ex]
&= (A+\mu V)^{-1}\frac{(A+ \mu V)(A+(\mu+h) V)^{-1}-I}{h}\\[1ex]
&=- (A+\mu V)^{-1}\,V\,(A+(\mu+h)V)^{-1}
\end{align*}
which tends to $- n_{\lambda}(\mu)\,y^\alpha \, n_{\lambda}(\mu)$ as $h\to 0$ in the  norm of $\mathcal B\left(L^p_m\right)$ since,  by Corollary \ref{Cor Rb},
$$ \mu \mapsto V(A+\mu)V)^{-1}=\mu y^\alpha \left (\mu-B^n +\frac{\lambda}{y^\alpha} \right )^{-1} \frac{1}{y^\alpha} $$
is continuous from $(0, \infty)$ to $\mathcal B\left(L^p_m\right)$. This shows  \eqref{Lem Mult Resolv eq 1} and then
$n_{\lambda}(\cdot)\in C^\infty\left(\R_+,\mathcal B\left(L^p_m\right)\right)$ and 
\begin{align}\label{Lem Mult Resolv eq 2}
D^k_\mu(n_{\lambda}(\mu))= a_k n_{\lambda}(\mu) \left(y^\alpha n_{\lambda}(\mu)\right)^{k},\qquad a_1=-1,\quad  a_{k+1}=-(k+1)a_k.
\end{align}
Formula \eqref{Lem Mult Resolv eq 2} follows  by induction after observing that, since $y^\alpha n_\lambda(\mu)$ and its derivative $D_\mu \left(y^\alpha n_\lambda\right)=- \left(y^\alpha \, n_{\lambda}(\mu)\right)^2$ commute, then $$D_\mu\left(y^\alpha n_\lambda(\mu)\right)^k=k D_\mu\left(y^\alpha n_\lambda(\mu)\right)\left(y^\alpha n_\lambda(\mu)\right) ^{k-1}=-k \left(y^\alpha n_\lambda(\mu)\right) ^{k+1}.$$

The  $\mathcal R$-boundedness of the family \eqref{Lem Mult Resolv eq 5} then follows from the  $\mathcal R$-boundedness of the families  \eqref{Lem Mult Resolv eq 3} since
\begin{align*}
	\mu^k D^k_\mu(\lambda n_{\lambda}(\mu))= a_k \lambda n_{\lambda}(\mu) \left(\mu y^\alpha n_{\lambda}(\mu)\right)^{k+1}.
\end{align*}
\qed

	In order to characterize  the domain of $y^{\alpha} B^n-y^{\alpha}$, we denote by $$D(y^\alpha)=\left\{u\in L^p_m:\ y^{\alpha}u\in L^p_m\right\}$$ the domain of the potential $V=y^\alpha$ in $L^p_m$.  Recalling  that Theorem \ref{Neumnan alpha m 2} assures that $D(y^{\alpha} B^n)=W^{2,p}_{\mathcal N}(\alpha,m)$, we consider, for $0<\frac{m+1}{p}<c+1-\alpha$, the Banach space
	\begin{align*}
		W^{2,p}_{\mathcal N}(\alpha,m)\cap D(y^\alpha)=\left\{u\in W^{2,p}_{loc}(\R_+):\ u,\ y^{\alpha}u,\    y^{\alpha}D_{yy}u,\ y^{\frac{\alpha}{2}}D_{y}u,\ y^{\alpha-1}D_{y}u\in L^p_m\right\}
	\end{align*}
endowed  with  norm $\|y^\alpha Bu\|_{L^p_m}+\|y^\alpha u\|_{L^p_m}+\| u\|_{L^p_m}$.

\begin{teo}  Let  $\alpha<2$, $\mu>0$, $c\in\R$.  Then for any $1<p<\infty$ such that $\alpha^-<\frac{m+1}p<c+1-\alpha$ the operator 	 
	 $L=y^{\alpha} B^n-\mu y^{\alpha}$ with domain $W^{2,p}_{\mathcal N}(\alpha,m)\cap D(y^\alpha)$ generates a bounded analytic semigroup  in $L^p_m$ which has maximal regularity. Moreover,
	
$$\mathcal {D}= \left\{u \in C_c^\infty ([0,\infty)): u \ {\rm constant \ in \ a \ neighborhood\  of }\  0  \right \}$$
is a core for $y^{\alpha} B^n-\mu y^{\alpha}$.
\end{teo}
{\sc{Proof.}}  The generation properties as well as the maximal regularity follows from Theorem \ref{Gen pert R}. Without any loss of generality we may assume that $\mu=1$.  We prove preliminarily that  $\mathcal {D}$ is dense in $W^{2,p}_{\mathcal N}(\alpha,m)\cap D(y^\alpha)=D(y^\alpha B^n) \cap D(y^\alpha)$. Let $u\in W^{2,p}_{\mathcal N}(\alpha,m)\cap D(y^{\alpha})$; up to  using a standard cut-off argument we may suppose that $\mbox{supp }u\subseteq [0,b]$ for some $b>0$. Using Remark \ref{Core compact}, let $(u_n)_{n\in\N}\subseteq  \mathcal D$ such that $\mbox{supp }u_n\subseteq [0,b]$ and  $u_n\to u$ in $W^{2,p}_{\mathcal N}(\alpha,m)$. Then by \cite[Proposition 3.2 (ii)]{MNS-Sobolev} 
\begin{align*}
	\|y^\alpha(u_n-u)\|_{L^p_m}\leq C \|y^{\alpha+1}(D_yu_n-D_yu)\|_{L^p_m}\leq  C b^2 \|y^{\alpha-1}D_y(u_n-u)\|_{L^p_m}
\end{align*}
which tends to $0$ as $n\to\infty$. This proves the density of $\mathcal {D}$.

Let us now characterize the domain. By definition  $D(y^\alpha B^n- y^{\alpha})=(1-y^\alpha B^n +y^{\alpha})^{-1}\left( L^p_m\right )$.  Let $u=(1-y^\alpha B^n +y^{\alpha})^{-1}f$ with $f \in L^p_m$. Using Corollary \ref{Cor Rb} (ii) we obtain
\begin{align}\label{eq spez}
	\|y^\alpha u\|_{L^p_{m}} +\|y^\alpha B u\|_{L^p_{m}}\leq C\left(\| \mathcal (y^\alpha B-y^\alpha) u\|_{L^p_{m}}+\| u\|_{L^p_{m}}\right)
\end{align}
which proves the inclusion $D(y^\alpha B^n-y^\alpha)\subseteq D(y^{\alpha} B^n)\cap D(y^{\alpha})$. To prove the reverse property we observe that, since the graph norm of $y^\alpha B^n-y^\alpha$ is clearly weaker than the norm of $D(y^{\alpha} B^n)\cap D(y^{\alpha})$,  inequality \eqref{eq spez} again shows that they are equivalent on 
$ D(y^\alpha B^n-y^\alpha)$, in particular on $ \mathcal D$ which is   dense in $D(y^{\alpha} B^n)\cap D(y^{\alpha})$, by the previous step. 
Therefore $ D(y^\alpha B^n-y^\alpha)=D(y^{\alpha} B^n)\cap D(y^{\alpha})$ and in particular $\mathcal D$ is a core.
\qed
We remark that Theorem \ref{Gen pert R} assures that $y^\alpha B^n-y^\alpha$  generates a semigroup on $L^p_m$  under the milder assumption  $0<\frac{m+1}p<c+1-\alpha$ and $c+1>0$. However, the hypothesis $(m+1)/p+\alpha>0$ must be added when $\alpha<0$ in order that $\mathcal D \subset D(y^\alpha)$.

\smallskip

The same method yields the domain of $B^n-\frac{1}{y^\alpha}$, using Corollary \ref{Cor Rb} (iii) with $m$ replaced by $m-\alpha p$ .

\begin{cor}
If $\alpha^+<\frac{m+1}{p} <c+1$, then the domain of $B^n-\frac{1}{y^\alpha}$ is $W^{2,p}_{\mathcal N}(0,m)\cap D(\frac{1}{y^\alpha})$.
\end{cor}

\bibliography{../TexBibliografiaUnica/References}

\end{document}